\documentclass{amsart}

\usepackage{color}
\definecolor{mygray}{gray}{.75}
\usepackage{xcolor}
\usepackage{paralist}
\usepackage{hyperref}
\usepackage{amssymb}
\usepackage{tikz,pgfplots}
\usepackage{bm}

\usepackage{pgfplotstable}

\pgfplotsset{compat=newest,compat/show suggested version=false}
\usepackage{cutwin}

\usepackage{dsfont}

\usepackage{algorithm}% http://ctan.org/pkg/algorithms
\usepackage{algpseudocode}% http://ctan.org/pkg/algorithmicx
\usepackage{caption}

\usepackage[font=footnotesize]{caption}

\makeatletter
\@namedef{subjclassname@2020}{%
  \textup{2020} Mathematics Subject Classification}
\makeatother

\theoremstyle{plain}
\newtheorem{thm}{Theorem}[section]
\newtheorem{rem}[thm]{Remark}
\newtheorem{prop}[thm]{Proposition}
\newtheorem{cor}[thm]{Corollary}

\newtheorem{example}[thm]{Example}

\newcommand{\cC}{\mathcal{C}}
\newcommand{\cD}{\mathcal{D}}
\newcommand{\cH}{\mathcal{H}}

\newcommand{\cL}{\mathcal{L}}
\newcommand{\cN}{\mathcal{N}}
\newcommand{\cP}{\mathcal{P}}
\newcommand{\cR}{\mathcal{R}}
\newcommand{\cT}{\mathcal{T}}
\newcommand{\cV}{\mathcal{V}}

\newcommand{\cX}{\mathcal{X}}

\newcommand{\bc}{\bm c} 
 
\newcommand{\bg}{\bm g} 
\newcommand{\br}{\bm r} 
 
\newcommand{\bu}{\bm u} 
\newcommand{\bv}{\bm v} 
\newcommand{\bw}{\bm w}
\newcommand{\bz}{\bm z} 
\newcommand{\bx}{\bm x} 
\newcommand{\by}{\bm y}

\newcommand{\bA}{\bm A} 
\newcommand{\bB}{\bm B}

\newcommand{\bG}{\bm G} 
\newcommand{\bI}{\bm I} 
\newcommand{\bM}{\bm M} 
\newcommand{\bN}{\bm N} 
\newcommand{\bQ}{\bm Q} 
\newcommand{\bR}{\bm R} 
\newcommand{\bU}{\bm U} 
\newcommand{\bV}{\bm V} 
\newcommand{\bW}{\bm W} 
\newcommand{\bX}{\bm X} 
\newcommand{\bY}{\bm Y} 
\newcommand{\bZ}{\bm Z}  
\newcommand{\bP}{\bm P}  

\newcommand{\bNull}{\mathbf{0}}

\newcommand{\N}{\mathbb{N}}
\newcommand{\R}{\mathbb{R}}
\newcommand{\U}{\mathbb{U}}
\newcommand{\V}{\mathbb{V}}

\newcommand{\bbB}{\pmb{\mathbb{B}}}
\newcommand{\bbK}{\pmb{\mathbb{K}}}
\newcommand{\bbM}{\pmb{\mathbb{M}}}

\newcommand{\Dt}{{\Delta t}}

\makeatletter
\newcommand{\opnorm}{\@ifstar\@opnorms\@opnorm}
\newcommand{\@opnorms}[1]{%
  \left|\mkern-1.5mu\left|\mkern-1.5mu\left|
   #1
  \right|\mkern-1.5mu\right|\mkern-1.5mu\right|
}
\newcommand{\@opnorm}[2][]{%
  \mathopen{#1|\mkern-1.5mu#1|\mkern-1.5mu#1|}
  #2
  \mathclose{#1|\mkern-1.5mu#1|\mkern-1.5mu#1|}
}
\makeatother
\makeatletter
\DeclareFontFamily{OMX}{MnSymbolE}{}
\DeclareSymbolFont{MnLargeSymbols}{OMX}{MnSymbolE}{m}{n}
\SetSymbolFont{MnLargeSymbols}{bold}{OMX}{MnSymbolE}{b}{n}
\DeclareFontShape{OMX}{MnSymbolE}{m}{n}{
    <-6>  MnSymbolE5
   <6-7>  MnSymbolE6
   <7-8>  MnSymbolE7
   <8-9>  MnSymbolE8
   <9-10> MnSymbolE9
  <10-12> MnSymbolE10
  <12->   MnSymbolE12
}{}
\DeclareFontShape{OMX}{MnSymbolE}{b}{n}{
    <-6>  MnSymbolE-Bold5
   <6-7>  MnSymbolE-Bold6
   <7-8>  MnSymbolE-Bold7
   <8-9>  MnSymbolE-Bold8
   <9-10> MnSymbolE-Bold9
  <10-12> MnSymbolE-Bold10
  <12->   MnSymbolE-Bold12
}{}

\let\llangle\@undefined
\let\rrangle\@undefined
\DeclareMathDelimiter{\llangle}{\mathopen}%
                     {MnLargeSymbols}{'164}{MnLargeSymbols}{'164}
\DeclareMathDelimiter{\rrangle}{\mathclose}%
                     {MnLargeSymbols}{'171}{MnLargeSymbols}{'171}
\makeatother
\usepackage{csquotes}

\DeclareMathOperator*{\Span}{span}

\numberwithin{equation}{section}

\pgfplotsset{select coords between index/.style 2 args={
    x filter/.code={
        \ifnum\coordindex<#1\fi
        \ifnum\coordindex>#2\fi
    }
}}

\usepackage{wrapfig}
\newsavebox{\wrapquestionfigure}
\opencutright

\title[Space-Time Form of the Wave Equation]{A Very Weak Space-Time Variational Formulation for the Wave Equation:\\ Analysis and Efficient Numerical Solution}
\author{Julian Henning}
\address{Ulm University, 
	Institute for Numerical Mathematics, 
	Helmholtzstr.\ 18, 89081 Ulm (Germany), 
	{\{julian.henning,karsten.urban\}@uni-ulm.de}}		
\author{Davide Palitta}
\address{Universit\`{a} di Bologna, 
	Centro AM$^2$, Dipartimento di Matematica, 
	Piazza di Porta S.\ Donato 5, 
	40127 Bologna (Italy),
	{\{davide.palitta,valeria.simoncini\}@unibo.it}}
\author{Valeria Simoncini} 
\author{Karsten Urban}

\thanks{The authors are grateful to Wolfgang Arendt (Ulm University) for various inspiring discussions concerning the analytical aspects. Furthermore, the authors acknowledge support by the state of Baden-Württemberg through bwHPC. 
The second and third authors are members of the Italian research group Indam-GNCS, whose support is gratefully acknowledged. 
Part of this work was carried out while the second author was affiliated with the Max Planck Institute for Dynamics of Complex
Technical Systems in Magdeburg, Germany.}
\date{Version of \today}

\begin{document}

\begin{abstract}
	We introduce a very weak space-time variational formulation for the wave equation, prove its well-posedness (even in the case of minimal regularity) and optimal inf-sup stability.  Then, we introduce a tensor product-style space-time Petrov-Galerkin discretization with optimal discrete inf-sup stability, obtained by a non-standard definition of the trial space. As a consequence, the numerical approximation error is equal to the residual, which is particularly useful for a posteriori error estimation.
	For the arising {discrete linear systems} in space and time, we  introduce efficient numerical solvers that appropriately exploit the equation structure, either at the preconditioning level or in the approximation phase by using a tailored Galerkin projection. This Galerkin method shows competitive behavior concerning {wall-clock} time, accuracy and memory as compared with a standard time-stepping method in particular in low regularity cases. Numerical experiments with a 3D (in space) wave equation illustrate our findings.
\end{abstract}

\maketitle

%===============================================
\section{Introduction}
\label{Sec:1}
%===============================================
The wave equation has extensively been studied in theory and numerical approximations. The aim of this paper is to introduce a (non-standard) variational Hilbert space setting for the wave equation {and a corresponding Petrov-Galerkin discretization} that is well-posed and optimally stable in the sense that the inf-sup constant is unity.  A major source of motivation for this view point is model reduction of parameterized partial differential equations by the reduced basis method, \cite{Haasdonk:RB,RozzaRB,QuarteroniRB}. {In that framework, the numerical approximation error is equal to the residual, which is particularly useful for a posteriori error estimation and model reduction.}

Space-time variational methods have been introduced, e.g., for parabolic problems \cite{r.andreev2013A,MR2891112,SpaceTimeUrbanPatera} and transport-dominated problems \cite{MR2084239,JBSU18,BTDeGh13,DHSW2012,DemGop11}, also partly with the focus of optimal inf-sup stability. The potential for efficient numerical solvers has been shown in \cite{HPSU20,Palitta2021}.

We follow the path of \cite{JBSU18,DHSW2012} and introduce a \emph{very weak} variational formulation in space and time by {applying} all derivatives onto the test functions using integration by parts. This means that the trial space is $L_2(I\times\Omega)$, where $I=(0,T)$ is the time interval and $\Omega\subset\R^d$ the domain in space. This is the \enquote{correct} space of minimal regularity for initial data $u_0\in L_2(\Omega)$.  Following \cite{DHSW2012}, we employ specifically chosen test spaces so as to derive  a well-posed variational problem. A Petrov-Galerkin method is then used for the discretization: inspired by \cite{JBSU18}, we first choose an appropriate test space and then define the (non-standard) trial space to preserve optimal inf-sup stability. 
This discretization results into a linear system of equations $\bbB_\delta \bu_\delta=\bg_\delta$, whose (stiffness) matrix $\bbB_\delta$ is a sum of tensor products and has large condition number, making the system solution particularly challenging. Memory and computational complexity are also an issue, as space-time discretizations in general lead to larger systems as compared to conventional time-stepping schemes, where a sequence of linear systems has to be solved, whose dimension corresponds to the spatial discretization only.

Building upon \cite{HPSU20}, we introduce matrix-based solvers that are competitive with respect to time-stepping schemes. In particular, we show that in case of minimal regularity the space-time method using fast matrix-based solvers outperforms a Crank-Nicolson time-stepping scheme.

The remainder of this paper is organized as follows: In Section \ref{Sec:2}, we review known facts concerning variational formulations in general and for the wave equation in particular. We derive an optimally inf-sup stable very weak variational form. Section \ref{Sec:3} is devoted to the Petrov-Galerkin discretization, again allowing for an inf-sup constant equal to 1. The arising linear system of equations is derived in Section \ref{Sec:4} and its efficient and stable numerical solution is discussed in Section \ref{Sec:5}. We show some results of numerical experiments for the 3D wave equation in Section~\ref{Sec:6}. For proving the well-posedness of the proposed variational form we need a result concerning a semi-variational formulation of the wave equation, whose proof is given in Appendix \ref{App:A}.

%===============================================
\section{Variational Formulations of the Wave Equation}
\label{Sec:2}
%===============================================
We are interested in a general linear equation of wave type. To this end, consider a Gelfand triple of Hilbert spaces $V\hookrightarrow H\hookrightarrow V'$ and a positive, symmetric operator $A\in\cL(D(A),H)$, where $D(A)$ is the domain of $A$ to be detailed in \eqref{eq:Domain} below.\footnote{We shall always denote by $V'$ the dual space of $V$ w.r.t.\ the pivot space $H$.} 
Setting $I:=(0,T)$, $T>0$ and given $g\in L_2(I; V')$\footnote{For a definition of Bochner spaces, see \S\ref{SubSec:Var} below.}, $u_0\in H$, $u_1\in V'$, we look for $u(t)\in V$, $t\in I\, \text{a.e.}$, such that
%--
\begin{equation}\label{Eq:1.1}
	\ddot{u}(t) + A\, u(t) = f(t)\,\, \text{in}\, V',\ t\in I\, \text{a.e.}, 
	\qquad 
	u(0)=u_0\in H,\, 
	\dot{u}(0)=u_1\in V'.
\end{equation}	
%--
Note, that the initial state is only in $H$ (e.g.\ $L_2(\Omega)$) and the initial velocity only in $V'$ (e.g.\ $H^{-1}(\Omega)$), which means very low regularity. Thus, without additional regularity, we cannot expect to get a smooth solution of \eqref{Eq:1.1}. Such non-smooth data are in fact a physically relevant situation. 
We restrict ourselves to LTI systems even though most of our results can be extended to the more general situation of a time-dependent operator $A(t)$. 

%----------------------------------------------------------------------------------
\subsection{Inf-sup-theory}\label{Sec:InfSup}
%----------------------------------------------------------------------------------
We are interested in finding a well-posed weak (or variational) formulation of \eqref{Eq:1.1}, i.e., Hilbert spaces $\U$, $\V$ of functions and a bilinear form $b:\U\times\V\to\R$ such that 
\begin{equation}\label{eq:VarFormGen}
	b(u,v) = g(v)\quad\forall v\in\V,
\end{equation}
has a unique solution $u\in\U$ for all given functionals $g\in\V'$ and that $u$ solves  \eqref{Eq:1.1} in some appropriate weak sense. The well-posedness of \eqref{eq:VarFormGen} is fully described by the following well-known fundamental statement. 

%-----------------------------------------------------------------
\begin{thm}[Ne\v{c}as Theorem, e.g.\ {\cite[Thm.\ 2]{NSV}}]\label{thm:Necas}
	Let $\U$, $\V$ be Hilbert spaces, let $g\in\V'$ be given and $b:\U\times\V\to\R$ be a bilinear form, which is bounded, i.e.
	\begin{align}
		\tag{C.1}
		&\exists\; \gamma<\infty:\quad
		b(u,v) \le \gamma \| u\|_\U\, \| v\|_\V,
		\quad\text{for all } u\in\U, v\in\V
		\quad
		\text{(boundedness).}
	\end{align}
	Then, the variational problem \eqref{eq:VarFormGen} admits a unique solution $u^*\in\U$, which depends continuously on the data $g\in\V'$ if and only if
	\begin{align}
		\tag{C.2}
		&\beta:=\inf_{u\in\U} \sup_{v\in\V} \frac{b(u,v)}{\| u\|_\U\, \| v\|_\V}>0 
		\quad\text{(inf-sup-condition)};
		\\
		\tag{C.3}
		&\forall\, 0\ne v\in\V\quad \exists\, u\in\U: \quad b(u,v)\ne 0
		\quad\text{(surjectivity)}.
	\end{align}
		\vskip-15pt\qed
\end{thm}
%-----------------------------------------------------------------

The inf-sup constant $\beta$ (or some lower bound) also plays a crucial role for the numerical approximation of the solution $u\in\U$ since it enters the relation of the approximation error and the residual (by the Xu-Zikatanov lemma \cite{MR1971217}, see also below). This motivates our interest in the size of $\beta$: the closer to unity, the better.

A standard tool (at least) for (i) proving the inf-sup-stability in (C.2); (ii) stabilizing finite-dimensional discretizations; and (iii) getting  sharp bounds for the inf-sup constant; is to determine the so-called \emph{supremizer}. To define it, let $b:\U\times\V\to\R$ be a generic bounded bilinear form and $0\ne u\in\U$ be given. Then, the \emph{supremizer} $s_u\in\V$ is defined as the unique solution of 
\begin{equation}\label{eq:supremizer}
	(s_u,v)_\V = b(u,v)\qquad \forall v\in\V.
\end{equation}
It is easily seen that
\begin{equation}\label{eq:supremizer2}
	\sup_{v\in\V} \frac{b(u,v)}{\| v\|_\V} = \sup_{v\in\V} \frac{(s_u,v)_\V}{\| v\|_\V} = \| s_u\|_\V,
\end{equation}
which justifies the name \emph{supremizer}. 

%----------------------------------------------------------------------------------
\subsection{The semi-variational framework}
\label{SubSec:SemiVar}
%----------------------------------------------------------------------------------
We start presenting some facts from the analysis of semi-variational formulations of the wave equation,  where we follow and slightly extend \cite[Ch.\ 8]{AU:eng}. The term \emph{semi-variational} originates from the use of classical differentiation w.r.t.\ time and a  variational formulation in the space variable. As above, we suppose that two real Hilbert spaces $V$ and $H$ are given, such that $V$ is compactly imbedded in $H$.  Let $a: V \times V \to \R$ be a continuous, coercive and symmetric bilinear form.\footnote{Note, that most of what is said can be also extended to $H$-elliptic forms (G\r{a}rding inequality).} 
Next, let $A$ be the operator on $H$ associated with $a(\cdot,\cdot)$ in the following sense: We define the \emph{domain} of $A$ by
\begin{align}\label{eq:Domain}
	D(A) := \{ u\in V:\, \exists f\in H \text{ such that } a(u,v) = (f,v)_H\, \forall v\in V\},
\end{align}
and recall that for any $u\in D(A)$ there is a unique $f\in H$ such that $a(u,v) = (f,v)_H$ for all $v\in V$. Then, we define
$A: D(A) \to H$ by $u\mapsto f:= Au$.
By the spectral theorem there exists an orthonormal basis $\{ e_n : n\in\N \}$ of $H$ and numbers $\lambda_n \in \R$ with $0< \lambda_1\leq \lambda_{2}\leq\cdots$, $\lim_{n\to\infty} \lambda_n=\infty$, such that 
\allowdisplaybreaks
\begin{subequations} 
	\begin{align}\label{eq:DA}
	V &= \Big\{ v\in H: \sum_{n=1}^{\infty} \lambda_n |(v,e_n)_H|^2 < \infty \Big\}, \\
	D(A) &=  \{ v\in H: Av\in H\}
	=\Big\{ v\in H: \sum_{n=1}^{\infty} \lambda_n^2 |(v,e_n)_H|^2 < \infty \Big\}, \\
	a(u,v) &= \sum_{n=1}^{\infty} \lambda_n (u,e_n)_H \, (e_n,v)_H, \qquad u,v\in V, \\
	A&v= \sum_{n=1}^{\infty} \lambda_n (v,e_n)_H \, e_n, \qquad v\in D(A). 
\end{align}
\end{subequations} 
In particular, $e_n \in D(A)$ and $A e_n = \lambda_n e_n$ for all $n \in \N$. For $s\in\R$, we define
\begin{align}\label{eq:Hs}
	H^s &:= \left\{ v = \sum_{n=1}^\infty v_n\, e_n:\, \| v\|_s^2 := \sum_{n=1}^\infty \lambda_n^s\, v_n^2<\infty\right\}
\end{align}
and note that $H^0=H$, $H^1=V$ and $H^2=D(A)$. Moreover, $(H^s)'\cong H^{-s}$, see Proposition \ref{Prop:A1}. 
We consider the non-homogeneous wave equation 
	\begin{align}\label{Eq:8.37neu}
		\ddot{w}(t) + A\, w(t) &= f(t), \quad t\in (0,T),  
		&&
		w(0) = u_0, \dot{w}(0) = u_1. 
	\end{align}
Then the following result on the existence and uniqueness holds. Its proof is given in Appendix \ref{App:A}.

\begin{thm}\label{Thm:A.2}
	Let $s\in\R_{\ge 0}$, $u_0 \in H^s$, $u_1 \in H^{s-1}$ and $f\in C([0,T]; H^{s-1})$. Then \eqref{Eq:8.37neu} admits a unique solution 
	\begin{align}\label{eq:cCs}
	w\in\cC^s:=C^2([0,T]; H^{s-2}) \cap C^1([0,T]; H^{s-1}) \cap C([0,T], H^s).
\end{align}
\end{thm}

We note a simple consequence for the backward wave equation.
\begin{cor}\label{Cor:A.4}
	Let $s\in\R_{\ge 0}$, $u_0 \in H^s$, $u_1 \in H^{s-1}$ and $g\in C([0,T]; H^{s-1})$. Then  
	\begin{align}\label{Eq:8.37rev}
		\ddot{w}(t) + A\, w(t) &= g(t), \quad t\in (0,T), 
	&&
		w(T) = u_0, \dot{w}(T) = u_1. 
	\end{align}
	admits a unique solution $w\in\cC^s$, see \eqref{eq:cCs}.
\end{cor}
\begin{proof}
	By the mapping $t\mapsto T-t$ we can transform \eqref{Eq:8.37rev} into \eqref{Eq:8.37neu} and deduce the well-posedness from Theorem \ref{Thm:A.2}.
\end{proof}

Theorem \ref{Thm:A.2} ensures that $B:= \frac{d^2}{dt^2} + A$ is an isomorphism of $\cC^s_0:=\{ v\in\cC^2: v(0)=\dot{v}(0)=0\}$ onto $C([0,T]; H^{s-2})$ for any $s\ge 0$. We detail the involved spaces in Table \ref{Tab:RegWave}, which also shows that we have to expect at most $w(t)\in H$, $t\in I$, in the semi-variational setting given the low regularity of the initial conditions in~\eqref{Eq:1.1}. Hence, in a variational space-time setting, we can only hope for $w(t)\in H$ for \emph{almost all} $t\in I$. 
%--------------------------------
\begin{table}[!hbt]
	\begin{center}
		\begin{tabular}{r|c|c|c||c|c|c}\hline
		$s$ & $u_0$ 	& $u_1$ & $f$ 	
				& $w$ 	& $\dot{w}$ 	& $\ddot{w}$ \\ \hline
		\footnotesize{$=$} & \multicolumn{2}{c|}{\footnotesize{$\in$}}
			& \footnotesize{$\in C([0,T];\cdot)$}
			& \multicolumn{3}{c}{\footnotesize{$\in C([0,T];\cdot)$}}
				\\ \hline\hline
		$0$ & $H$ 	& $V'$ 	& $V'$ 	& $H$ 	& $V'$ 		& $D(A)'$ \\ \hline
		$1$ & $V$ 	& $H$ 	& $H$ 	& $V$ 	& $H$ 		& $V'$ \\ \hline
		$2$ & $D(A)$ 	& $V$ 	& $V$ 	& $D(A)$ 	& $V$ 		& $H$ \\ \hline
		\end{tabular}
		\caption{\label{Tab:RegWave}Regularity statements for the wave equation -- classical in time, variational in space.}
	\end{center}
\end{table}

%----------------------------------------------------------------------------------
\subsection{Biharmonic problem and mixed form}
\label{Sec:Biharm}
%----------------------------------------------------------------------------------
For later reference, let us consider the bilinear form $q: D(A)\times D(A)\to\R$ defined by $q(u,v):= (Au,Av)_H$, $u,v\in D(A)$, which is of biharmonic type. {In order to detail the associated operator $Q$, recall that we have a Gelfand quintuple $D(A) \hookrightarrow V \hookrightarrow H \hookrightarrow V' \hookrightarrow D(A)'$.} The duality pairing of $D(A)$ and $D(A)'$ is denoted by $\langle\cdot,\cdot\rangle_{D(A)'\times D(A)}$. Then, $Q:D(A)\to D(A)'$ defined as $\langle Q u,v \rangle_{D(A)'\times D(A)}=q(u,v)$ for $u,v\in D(A)$.

The adjoint operator $A':H\to D(A)'$ is given by $\langle A'h,w\rangle_{D(A)'\times D(A)}= (h,Aw)_H$ for $w\in D(A)$ and $h\in H$. Then, $A'A: D(A)\to D(A)'$ and we get for $u,v\in D(A)$ that
$\langle A'A u,v \rangle_{D(A)'\times D(A)}
	= (Au, Av)_H 
	= q(u,v) 
	= \langle Q u,v \rangle_{D(A)'\times D(A)}$,
hence $Q=A'A$. Next, we consider the following operator problem: 
\begin{align}\label{eq:Biharm}
	\text{given } 
		g\in D(A)', 
	\text{ determine }
	z\in D(A) 
	\text{ such that }
	Qz=g. 
\end{align}
Introducing the auxiliary variable $u:=Az\in H$, we can rewrite this problem as
\begin{align}\label{eq:SPPBiharm}
	\begin{pmatrix}
		I & A \\ A' & 0 
	\end{pmatrix}
	\begin{pmatrix} u \\ z \end{pmatrix}
	= 
	\begin{pmatrix} 0 \\ -g \end{pmatrix},
\end{align}
which is easily seen to be equivalent to \eqref{eq:Biharm}. 

%--------------------------------------------------------------------------------
\subsection{Towards space-time variational formulations}
\label{SubSec:Var}
%--------------------------------------------------------------------------------
The semi-variational formulation described above cannot be written as a variational formulation in the form of \eqref{eq:VarFormGen}, since $C^k([0,T]; X)$ is not a Hilbert space, even if  $X$ is a Hilbert space of functions $\phi:\Omega\to\R$ in space, e.g.\ $L_2(\Omega)$ or $H^1_0(\Omega)$. We need Lebesgue-type spaces for the temporal and spatial variables yielding the notion of \emph{Bochner spaces}, denoted by $\cX:=L_2(I; X)$\footnote{Spaces of space-time functions are denoted by calligraphic letters, spaces of functions in space only by plain letters.} and defined as
$$
	\cX 
	:= L_2(I; X) 
	:= \bigg\{ v: I\to X:\, \| v\|_{L_2(I; X)}^2 := \int_0^T \| v(t)\|_X^2\, dt <\infty\bigg\},
$$
which are Hilbert spaces with the inner product 
$(w,v)_\cX := \int_0^T ( w(t), v(t))_X\, dt$,
where $(\cdot,\cdot)_X$ denotes the respective inner product in $X$. We will often use the specific cases $(\cdot,\cdot)_\cV$ and $(\cdot,\cdot)_\cH$ for $\cV:=L_2(I; V)$ as well as $\cH:=L_2(I;H)$. Sobolev-Bochner spaces, e.g.\ $H^1(I;X)$, $H^2(I;X)$ can be defined accordingly using weak derivatives w.r.t.\ the time variable.

We will derive a space-time variational formulation in Bochner spaces, i.e., we multiply the partial differential equation  in \eqref{Eq:1.1} with test functions in space \emph{and} time and also integrate w.r.t.\ both variables. Now, the question remains how to apply integration by parts. One could think of performing integration by parts once w.r.t.\ all variables. This would yield a variational form in the Bochner space $H^1(I; V)$. However, we were not able to prove well-posedness in that setting.  Hence, we suggest a \emph{very} or \emph{ultra weak} variational form, where all derivatives are put onto the test space by means of integration by parts. We thus define the trial space as
\begin{align}\label{eq:defW}
	\U := \cH = L_2(I;H)
\end{align}
and search for an appropriate test space $\V$ to guarantee well-posedness of~\eqref{eq:VarFormGen}. Performing integration by parts twice for both time and space variables, we obtain
	\allowdisplaybreaks[1]
	\begin{align}
		\label{eq:veryweak}
		b(u,v) := (u, \ddot{v} + Av)_\cH, 
	&&
		g(v)	:= (f,v)_\cH 
					+ \langle u_1, v(0)\rangle
					- (u_0, \dot{v}(0))_H, 	
	\end{align}
for $v\in\V$, where the space $\V$ still needs to be defined in such a way that all assumptions of Theorem \ref{thm:Necas} are satisfied. It turns out that this is not a straightforward task. The duality pairing $\langle\cdot,\cdot\rangle$ is defined in \eqref{eq:Duality} in the appendix.

\subsubsection*{The Lions-Magenes theory}
Variational space-time problems for the wave equation within the setting \eqref{eq:veryweak} have already been investigated in the book \cite{LM1} by Lions and Magenes. We are going to review some facts from  \cite[Ch.\ III, \S9, pp.\ 283-299]{LM1}. The point of departure is the following adjoint-type problem.

For a given $\varphi\in L_2(I; H)=\U$, find $v:I\times\Omega\to\R$ such that
\begin{equation}\label{LM:9.1}
	\ddot{v} + A\, v = \varphi,
	\qquad v(T)=\dot{v}(T)=0.
\end{equation}
It has been shown that the following space\footnote{The definition \eqref{LM:9.2} is literally cited from \cite{LM1}.}
\begin{equation}\label{LM:9.2}
	\V := \text{space described by the solution $v$ of \eqref{LM:9.1} as $\varphi$ describes $L_2(I;H)$}
\end{equation}
plays an important role for the analysis. 
It is known that $\V \subset C([0,T]; V) \cap C^1([0,T]; H) \cap H^2(I;V')$ and that $\frac{d^2}{dt^2}+A$ is an isomorphism of $\V$ onto $\U$. 

%-----------------------------------------------------------------
\begin{thm}\label{Thm:LM-1} \cite[Ch.\ 3, Thm.\ 8.1, 9.1]{LM1}
	Let $a:V\times V\to\R$ satisfy a G\r{a}rding inequality and let $f\in L_2(I;H)$, $u_0\in V$, $u_1\in H$ be given. Then, 
	\begin{compactenum}
	\item[(a)] 
	there is a unique $u^*\in H^1(I;H)\cap L_2(I;V)$ such that $\ddot{u}^* + Au^*=f$, $u^*(0)=u_1$, $\dot{u}^*(0)=u_1$. 
	In addition $u^*\in H^2(I;V')$;
	\item[(b)] for any $\ell\in\V'$ there is a unique $u^*\in\U$ such that $b(u^*,v) = \ell(v)$ for all $v\in\V$.\hfill\qed
	\end{compactenum}
\end{thm}
%-----------------------------------------------------------------
Notice that the first statement is proven by deriving energy-type estimates for the uniqueness and a Faedo-Galerkin approximation for the existence. Let us comment on the previous theorem. First, we note that $u_0\in V$, $u_1\in H$ are `too smooth' initial conditions, we aim at (only) $u_0\in H$, $u_1\in V'$, see \eqref{eq:veryweak}. As a consequence:
	\begin{compactenum}
		\item Statement (a) in Thm.\ \ref{Thm:LM-1} results in a `too smooth' solution. In fact, we are interested in a very weak solution $u\in L_2(I;H)$, (a) is `too' much.
		\item Even though the stated solution in (b) has the `right' regularity, it is not clear how to associate the functional $g$ in \eqref{eq:veryweak} to the dual space $\V'$, i.e., how to interpret the three terms of $g$ in \eqref{eq:veryweak} in the space $\V'$.
	\end{compactenum}
These issues are partly fixed by the following statement.

%-----------------------------------------------------------------
\begin{thm}\label{Thm:LM-2} \cite[Ch.\ III, Thm.\ 9.3, 9.4]{LM1}
	Let the bilinear form $a(\cdot,\cdot)$ be coercive, $f\in L_2(I;V')$, $u_0\in H$, $u_1\in V'$. Then, there exists a unique $u^*\in L_\infty(I;H)\cap W^{1}_\infty(I;V')$ such that $b(u^*,v) = g(v)$ for all $v\in \V_0 := \V\cap L_2(I;W)$ 
	with $b(\cdot,\cdot)$ and $g$ defined as in \eqref{eq:veryweak}. Moreover, $u^*\in C^0(\bar{I};H)\cap C^1(\bar{I};V')$.\hfill\qed
\end{thm}
%-----------------------------------------------------------------

Even though the latter result uses the `right' smoothness of the data and also includes existence \emph{and} uniqueness, we are not fully satisfied with regard to our goal of a well-posed variational formulation of the wave equation in Hilbert spaces. In fact, the `trial space' $L_\infty(I;H)\cap W^{1}_\infty(I;V')$ is not a Hilbert space and it is at least not straightforward to see how we can base a Petrov-Galerkin approximation on such a trial space. Hence, we follow a different path.

%----------------------------------------------------------------------------------
\subsection{An optimally inf-sup stable very weak variational form}
%----------------------------------------------------------------------------------
We are going to derive a well-posed very weak variational formulation \eqref{eq:VarFormGen} of \eqref{Eq:1.1}, where $\U=L_2(I;H)$ and $b(\cdot,\cdot)$, $g(\cdot)$ are defined by \eqref{eq:veryweak}. To this end, we will follow the framework presented in \cite{DHSW2012}. This approach is also called the \emph{method of transposition} and -- also -- already goes back to \cite{LM1}, see also e.g.\ \cite{MR3535626,MR2084239,MR3070527} for the corresponding finite element error analysis. For the presentation we will need the semi-variational formulation described above. 

Let us restrict ourselves to $A=-\Delta$ acting on a convex domain $\Omega\subset\R^d$ and supplemented by homogeneous Dirichlet boundary conditions. This means that $H=L_2(\Omega)$, $V=H^1_0(\Omega)$ and $D(A)=H^2(\Omega)\cap H^1_0(\Omega)$. However, we stress the fact that most of what is said here can be also extended to other elliptic operators.
Then, the starting point is the operator equation in the classical form, i.e.,
$$
	B_\circ u = g,
	\quad\text{where }\,
	B_\circ = \frac{d^2}{dt^2}+A_\circ,
	\qquad
	\Omega_T := (0,T) \times \Omega, 
$$
i.e., $A_\circ=-\Delta$ is also to be understood in the classical sense. Next, denote the classical domain of $B_\circ$ by $\cD(B_\circ)$, where initial and boundary conditions are also imposed in $\cD(B_\circ)$, i.e.,  $\cD(B_\circ):=\{ v\in C(\bar\Omega_T): B_\circ v \in C(\Omega_T), v(0)=0, v(t,\cdot)_{|\partial\Omega}=0\,\, \forall t\in [0,T]\}$. Hence,
$$
	\cD(B_\circ) 
		= C^2(\Omega_T) \cap C^1_{\{ 0\}} ([0,T]; C_0(\overline\Omega))
		\kern-1pt
		= \big[ C^2(I) \cap C^1_{\{ 0\}}([0,1])\big] \times \big[ C^2(\Omega)\cap C_0(\overline\Omega)\big],
$$
where $C_0(\overline\Omega) := \{ \phi\in C(\overline\Omega):\, \phi_{|\partial\Omega}=0\}$ models the homogeneous Dirichlet conditions, and for $t\in [0,T]$ and any pair of function spaces $X$, $Y$, we define 
\begin{align*}
	C^1_{\{ t\}} ([0,T]; X,Y) := \{ u\in C([0,T]; X)\cap C^1([0,T]; Y) :\, u(t)=0,\, \dot{u}(t)=0\}.
\end{align*}
The range $\cR(B_\circ)$ in the classical sense then reads $\cR(B_\circ) = C(\overline\Omega_T)$. 
As a next step, we determine the formal adjoint $B_\circ^*$ of $B_\circ$. Since
\begin{align*}
	(B_\circ u, v)_\cH = (u, B_\circ v)_\cH
	\quad\text{ for all } u,v\in C^\infty_0(\Omega_T),
\end{align*}
the operator $B_\circ$ is self-adjoint -- but  with homogeneous terminal conditions $u(T)=\dot{w}(T)=0$ instead of initial conditions. This means that $\cR(B_\circ^*) = C(\overline\Omega_T)$ and
\begin{align*}
	\cD(B_\circ^*) 
		&= C^2(\Omega_T) \cap C^1_{\{ T\}} ([0,T]; C_0(\overline\Omega)) 
		\kern-3pt=\kern-3pt \big[ C^2(I) \cap C^1_{\{ T\}}([0,1])\big] \times \big[ C^2(\Omega)\cap C_0(\overline\Omega)\big].
\end{align*}
Following \cite{DHSW2012}, we need to verify the following conditions
\begin{compactenum}
	\item[($B^*1$)] $B_\circ^*$ is injective on the dense subspace $\cD(B_\circ^*)\subset L_2(I;H)$ and 
	\item[($B^*2$)] $\cR(B_\circ^*)\hookrightarrow L_2(I;H)$ is densely imbedded.
\end{compactenum}
Since $C(\overline\Omega_T)\cong C([0,T]; C(\overline\Omega)) \hookrightarrow L_2(I;H)$ is dense, ($B^*2$) is immediate. In order to prove ($B^*1$), first note that 
\begin{align}\label{eq:DefC2T}
	\cD(B_\circ^*)\subset \cC^2_{T} 	
	:= \cC^2 \cap C^1_{\{T\}}([0,T]; V).
\end{align}
Let us denote the continuous extension of $B_\circ^*$ from $\cD(B_\circ^*)$ to $\cC^2_{T}$ also by $B_\circ^*$. Corollary \ref{Cor:A.4} implies that this continuous extension $B_\circ^*$ is an isomorphism from $\cC^2_{T}$ onto $C([0,T]; V)$ (here we need the semi-variational theory). This implies that $B_\circ^*$ is injective on $\cD(B_\circ^*)$, i.e., ($B^*1$). 
Now, the properties ($B^*1$) and ($B^*2$) ensure that 
\begin{align}\label{eq:normV}
	\| v\|_\V := \| B_\circ^* v\|_\cH
\end{align}
is a norm on $D(B_\circ^*) = \cC_{T}^2$. Then, we set
\begin{align}\label{Def:V}
	\V := \mathrm{clos}_{\|\cdot\|_\V} (\cC^2_{T}) \subset  L_2(I;H),
	\quad
	(v,w)_\V := (B^* v, B^*w)_\cH,
	\, v,w\in\V,
\end{align}
which is a Hilbert space, where $B^*$ is to be understood as the continuous extension of $B^*_\circ$ from $\cC^2_{T}$ to $\V$.  
Now, we are ready to prove our first main result.

%------------------------------------------------------------------------------------------------------------------------------
\begin{thm}\label{Thm:stablecont}
	Let $f\in L_2(I;V')$, $u_0\in H$ and $u_1\in V'$. Moreover, let $\V$, $b(\cdot,\cdot)$, and $g(\cdot)$ be defined as in \eqref{Def:V} and \eqref{eq:veryweak}, respectively. 
	Then, the variational problem
	\begin{equation}\label{Eq:vwEqt}
		b(u,v) = g(v) 
		\quad\text{ for all } v\in\V,
	\end{equation}
	admits a unique solution $u^*\in\U$. 
	In particular,
	\begin{equation}
		\beta := 
		\inf_{u\in\U} \sup_{ v\in\V} \frac{b(u, v)}{\| u\|_\U\, \|  v\|_{\V}} 
		= \sup_{u\in\U} \sup_{ v\in\V} \frac{b(u, v)}{\| u\|_\U\, \|  v\|_{\V}} 
		= 1.
	\end{equation}
\end{thm}
%------------------------------------------------------------------------------------------------------------------------------
\begin{proof}
	We are going to show the conditions (C.1)-(C.3) of Theorem \ref{thm:Necas} above.\\
	(C.1) \emph{Boundedness}: Let $u\in\U$, $ v\in\V$, then by Cauchy-Schwarz' inequality
	\begin{align*}
		b(u, v)
		&= (u,\ddot{ v}+A v)_\cH 
		\le \| u\|_\cH \, \| \ddot{ v} + A v\|_\cH
		= \|u\|_\U\, \| v\|_{\V},
	\end{align*}
	i.e., the continuity constant is unity. \\
	(C.2) \emph{Inf-sup}: Let $0\ne u\in\U$ be given. We consider the supremizer $s_u\in \V$ defined as $(s_u,  v)_{\V} = b(u, v) = (u,\ddot{ v}+A v)_\cH$ for all $ v\in \V$. Since  by definition of the inner product 
	$
	(s_u,  v)_{\V} 
	= ( \ddot{s}_u, \ddot{v})_\cH + (As_u, A v)_\cH 
	$ 
	for all $ v\in \V$ we get $\ddot{s}_u + As_u=u$ in $\cH$. Then, by \eqref{eq:supremizer2}, 
	\begin{align*}
	\sup_{ v\in\V} \frac{b(u, v)}{\|  v\|_{\V}} 
	&= \sup_{ v\in\V} \frac{(s_u, v)_{\V}}{\|  v\|_{\V}} 
	= \| s_u\|_{\V}
	= \| \ddot{s}_u + As_u\|_\cH 
	= \| u\|_\cH,
	\end{align*}
	i.e., $\beta=1$ for the inf-sup constant. \\
	(C.3) \emph{Surjecitivity}: Let $0\ne v\in\V$ be given. Then, there is a sequence $(v_n)_{n\in\N}\subset\cC_{T}^2$ with $v_n\not=0$, converging towards $v$ in $\V$. Since $B_\circ^*$ is an isometric isomorphism of $\cC^2_{T}$ onto $C([0,T];V)$, there is a unique $u_n:=B^*_\circ v_n=\ddot{v}_n + Av_n\in C([0,T];V)$. Hence $0\not= \| v_n\|_{\cC^2} = \| u_n\|_{C([0,T];V)}$. Possibly by taking a subsequence,  $(u_n)_{n\in\N}$ converges to a unique limit $u_v\in L_2(I;H)$. 
	We take the limit as $n\to\infty$ on both sides of $u_n=\ddot{v}_n + Av_n$ and obtain $0\ne u_v=B^*v=\ddot{v} + Av\in L_2(I;H)=\U$. 
	Finally, $b(u_v,v) = (u_v, B^*v)_\cH = (u_v, u_v)_\cH = \| u_v\|_\U^2> 0$, which proves surjectivity and concludes the proof.
\end{proof}

\begin{rem}
	The essence of the above proof is the fact that $\U$ and $\V$ are related as $\U=B^*(\V)$ and noting that $B$ is self-adjoint up to initial versus terminal conditions.
\end{rem}

\subsubsection*{Further remarks on the test space $\V$}
The above definition \eqref{Def:V} is not well suited for a discretization. Hence, we are now going to further investigate $\V$. First, note that $\V = \mathrm{clos}_{\|\cdot\|_\V} (\cD(B_\circ^*))$ and recall that $\cD(B_\circ^*)= \big[ C^2(I) \cap C^1_{\{ T\}}([0,1])\big] \times \big[ C^2(\Omega)\cap C_0(\overline\Omega)\big]$ is a tensor product space. Next, for $v(t,x) = \vartheta(t)\, \varphi(x)$, we define a tensor product-type norm as
\begin{align*}
	\opnorm{v}_\V^2 
	= \opnorm{\vartheta\otimes\varphi}_\V^2 
	&:= \big( \|\ddot\vartheta\|_{L_2(I)}^2 + \|\vartheta\|_{L_2(I)}^2\big)
		\, \big( \|\varphi\|_{H}^2 + \|A\varphi\|_{H}^2\big) 
	=: \opnorm{\vartheta}_t^2\, \opnorm{\varphi}_x^2,  
\end{align*}
and set (for $A=-\Delta$ on $\Omega\subset\R^d$)
\begin{align}
	\V_\circ
	&:= \mathrm{clos}_{\opnorm{\cdot}_\V}  (\cD(B_\circ^*)) 
		\nonumber\\
	&= \mathrm{clos}_{\opnorm{\cdot}_t}  \big(C^2(I) \cap C^1_{\{ T\}}([0,1])\big)
		\times \mathrm{clos}_{\opnorm{\cdot}_x}  \big(C^2(\Omega)\cap C_0(\overline\Omega)\big) \nonumber \\
	&= H^2_{\{ T\}}(I) \times \big[ H^2(\Omega) \times H^1_0(\Omega)\big], \label{eq:Vcirc}
\end{align}
where  $H^2_{\{T\}}(I) := \{ \vartheta\in H^2(I)\!:\, \vartheta(T)=\dot\vartheta(T)=0\}$ recalling that $D(-\Delta)=H^2(\Omega) \times H^1_0(\Omega)$, \cite{MR2084239}. Again, it is readily seen that $\V_\circ \subset \V$, but the contrary is not true in general. In view of \eqref{eq:Vcirc}, $\V_\circ$ is a tensor product space which can be discretized in a straightforward manner.

%===============================================
\section{Petrov-Galerkin Discretization}
\label{Sec:3}
%===============================================
We determine a numerical approximation to the solution of a variational problem of the general form \eqref{eq:VarFormGen}.  To this end, one chooses finite-dimensional trial and test spaces, $\U_\delta\subset\U$, $\V_\delta\subset\V$, respectively, where $\delta$ is a discretization parameter to be explained later.  For convenience, we assume that their dimension is equal, i.e., $\cN_\delta:=\dim\U_\delta=\dim\V_\delta$. The Petrov-Galerkin method then reads
\begin{align}\label{eq:var-disc}
	\text{find } u_\delta\in\U_\delta:\quad
	b(u_\delta,v_\delta) =   g(v_\delta)
	\quad\text{for all } v_\delta\in \V_\delta.
\end{align}
As opposed to the coercive case, the well-posedness of \eqref{eq:var-disc} is not inherited from that of \eqref{Eq:vwEqt}. In fact, in order to ensure uniform stability (i.e., stability independent of the discretization parameter $\delta$), the spaces $\U_\delta$ and $\V_\delta$ need to be appropriately chosen in the sense that the discrete inf-sup (or LBB -- Ladyshenskaja-Babu\v{s}ka-Brezzi) condition holds, i.e., there exists a $\beta_\circ>0$ such that
\begin{align}\label{eq:LBB}
	\beta_\delta
	&:=\inf_{u_\delta\in\U_\delta} \sup_{v_\delta\in\V_\delta} 
		\frac{b(u_\delta,v_\delta)}{\| u_\delta\|_\U\, \| v_\delta\|_\V}
	\ge \beta_\circ >0,
\end{align}
where the crucial point is that $\beta_\circ$ is independent of $\delta$. The size of $\beta_\circ$ is also relevant for the error analysis, since the Xu-Zikatanov lemma \cite{MR1971217} yields a best approximation result 
\begin{align}\label{eq:bestApprox}
	\| u^*-u^*_\delta\|_\U \le \frac{1}{\beta_\circ} \inf_{w_\delta\in\U_\delta}\| u^*-w_\delta\|_\U
\end{align}
for the `exact' solution $u^*$ of \eqref{Eq:vwEqt} and the `discrete'  solution $u_\delta^*$ of \eqref{eq:var-disc}. This is also the key for an optimal error/residual relation, which is important for a posteriori error analysis (also within the reduced basis method).

%---------------------------------------------------------------------------------------------------------
\subsection{A stable Petrov-Galerkin space-time discretization}
%---------------------------------------------------------------------------------------------------------
To properly discretize $\V$,
we consider the tensor product subspace $\V_\circ \subset \V$ introduced in \eqref{eq:Vcirc} which allows for a straightforward finite element discretization. Hence, we look for a pair $\U_\delta\subset\U$ and $\V_\delta\subset\V_\circ$ satisfying \eqref{eq:LBB} with a possibly large inf-sup lower bound $\beta_\circ$, i.e., close to unity. Constructing such a stable pair of trial and test spaces is again a nontrivial task, not only for the wave equation. It is a common approach to choose some trial approximation space $\U_\delta$ (e.g.\ by splines) and then (try to) construct an appropriate according test space $\V_\delta$ in such a way that \eqref{eq:LBB} is satisfied. This can be done, e.g., by computing the supremizers for all basis functions in $\U_\delta$ and then define $\V_\delta$ as the linear span of these supremizers. However, this would amount to solve the original problem $\cN_\delta$ times, which is way too costly. We mention that this approach indeed works within the discontinuous Galerkin (dG) method, see, e.g., \cite{BTDeGh13,DemGop11}.
We will follow a different path, also used in \cite{JBSU18} for transport problems. We first construct a test space $\V_\delta$ by a standard approach and then define a stable trial  space $\U_\delta$ in a second step. This implies that the trial functions are no longer `simple' splines but they arise from the application of the adjoint operator $B^*$ (which is here the same as the primal one $B$ except for initial/terminal conditions) to the test basis functions.

%%%%%%%%%%%%%%%%%%%%%%%%%%%%%%%
\subsubsection*{Finite elements in time.} 
We start with the temporal discretization. 
We choose some integer ${N_t}>1$ and set $\Dt:=T/{N_t}$. This results in a temporal \enquote{triangulation} 
\begin{align*}
	\cT_{\Dt}^\text{time}\equiv\{ t^{k-1}\equiv(k-1)\Delta t < t \le k\, \Dt \equiv t^k, 1\le k\le {N_t}\}
\end{align*}
in time. Then, we set 
\begin{align}\label{eq:RDt}
	R_\Dt := \Span\{ \varrho^1,\ldots ,\varrho^{{N_t}}\}\subset H^2_{\{T\}}(I), 
\end{align}	
e.g.\ piecewise quadratic splines on $\cT_{\Dt}^\text{time}$ with standard modification in terms of multiple knots at the right end point of $\bar I=[0,T]$.

\begin{example}\label{Ex:Time}
Denote by $S^k$ the quadratic B-spline corresponding to the nodes $t^{k-2}$, $t^{k-1}$, $t^{k}$ and $t^{k+1}$, where we extend the node sequence outside $\overline{I}$ in an obvious manner. Then, $\varrho^{k}:=S^{k-1}$, $k=3,...,N_t$ are $H^2_0(I)$-functions which are fully supported in $I$. The remaining two basis functions on the left end point of the interval $I$, i.e., $\varrho^{1}$, $\varrho^2$, can be formed by using $t^0=0$ as double and triple node, respectively. Thus, we get a discretization in $H^2_{\{T\}}(I)$ of dimension $N_t$. We show an example for $T=1$ and $\Dt=\frac18$ (i.e., $N_t=8$) in Figure \ref{1D:Examples}, the test functions in the center, optimal trial functions on the right.
\end{example}

\subsubsection*{Discretization in space.} 
For the space discretization, we  choose any conformal finite element space
\begin{align}\label{eq:Zh}
	Z_h :=\Span\{ \phi_1,\ldots,\phi_{{N_h}}\}\subset H^1_0(\Omega)\cap H^2(\Omega),
\end{align}
e.g.\ piecewise quadratic finite elements with homogeneous Dirichlet boundary conditions.

\begin{figure}[!tb]
	\begin{center}
			\includegraphics[width=0.32\textwidth]{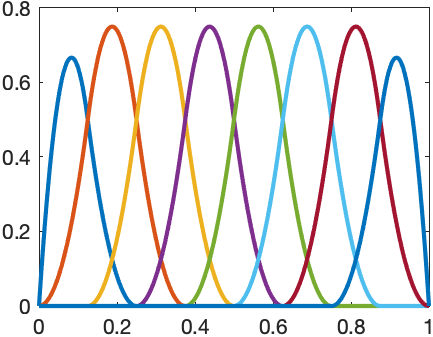}
			\includegraphics[width=0.32\textwidth]{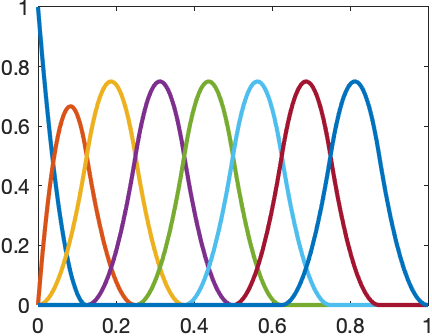}
			\includegraphics[width=0.32\textwidth]{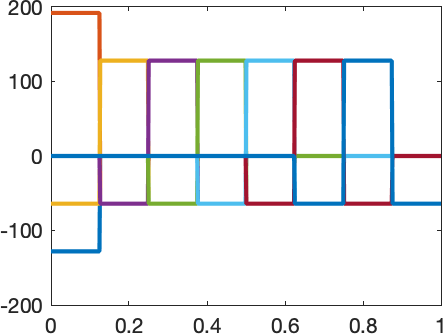}
		\caption{\label{1D:Examples}Discretization for time and 1D-space, $h=\Dt=\frac18$. 
		Quadratic B-spline test functions, from left to right: 1D-space, time and inf-sup-optimal trial functions.}
	\end{center}
\end{figure}

\begin{example}\label{Ex:Space1d}
	As an example for the space discretization, let us detail the univariate (1D) case $\Omega=(0,1)$. Define $x_j:=j\, h$, $j=0,...,N_h:=\frac1h$, and denote by $S^j$ the quadratic B-spline corresponding to the nodes $x_{j-2}$, $x_{j-1}$, $x_{j}$, and $x_{j+1}$. The B-splines $\phi_j:=S_{j+1}$, $j=2,...,N_h-1$, are supported in $\Omega$. We define the two boundary functions $\phi_1$ and $\phi_{N_h}$ as  the quadratic B-spline w.r.t.\ the nodes $(0,0,x_1,x_2)$ and $(x_{N_h-2},x_{N_h-1},1,1)$ (i.e., with double nodes), respectively, such that the homogeneous boundary conditions are satisfied. 
	We obtain a discretization of dimension $N_h$. We show an example for $\Omega=(0,1)$ and $h=\frac18$ (i.e., $N_h=8$) in Figure \ref{1D:Examples}, the test functions on the left. The arising trial functions are depicted on the right and turn out to be identical with to time discretization in Example \ref{Ex:Time}.
\end{example}

\subsubsection*{Test and trial space in space and time.} 
Then, we define the test space as
\begin{align}\label{eq:Vdelta}
	\V_\delta &:= R_\Dt\otimes Z_h \subset \V_\circ \subset \V, 
	\qquad \delta=(\Dt,h), \\
	&= \Span\{ \varphi_\nu := \varrho^k\otimes \phi_i:\, k=1,...,N_t,\, i=1,...,N_h, \nu=(k,i)\}, 
	\nonumber 
\end{align}
which is a tensor product space of dimension $\cN_\delta = N_t\, N_h$.
\vspace*{5pt}

\hspace*{-\parindent}%
\begin{minipage}{0.54\textwidth}
The trial space $\U_\delta$ is constructed by applying the adjoint operator $B^*$ to each test basis function, i.e., for $\mu=(\ell,j)$ and $A=-\Delta$
\begin{align*}
	\psi_\mu 
	&:= B^*(\varphi_\mu) 
	=  B^*(\varrho^\ell\otimes \phi_j) \\
	&= (\partial_{tt} + A)(\varrho^\ell\otimes \phi_j)\\
	&= \ddot{\varrho}^\ell \otimes \phi_j + \varrho^\ell\otimes  A\phi_j,
\end{align*}
i.e., $\U_\delta := B^*(\V_\delta) = \Span\{ \psi_\mu:\, \nu=1,...,\cN_\delta\}$. Since $B^*$ is an isomorphism of $\V$ onto $L_2(I;H)$, the functions $\psi_\nu$ are in fact linearly independent. An example of a single trial function is shown in Figure \ref{Fig:Trial}.
\end{minipage}
\begin{minipage}{0.45\textwidth}
	\centering
		\includegraphics[width=.99\textwidth]{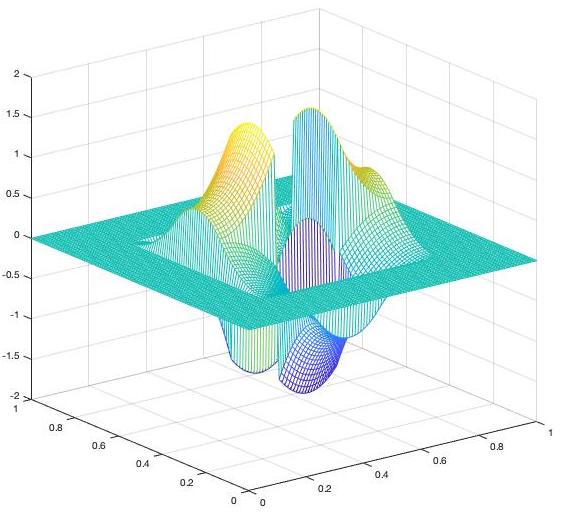}\vspace*{-15pt}
		\captionof{figure}{Sample trial function for $I=\Omega=(0,1)$, $\Dt=h=0.2$.}
		\label{Fig:Trial}
\end{minipage}
\vspace*{5pt}

\begin{prop}\label{Lem:LBB}
 	For the space $\V_\delta$ defined in \eqref{eq:Vdelta} and $\U_\delta:=B^*(\V_\delta)$, we have
	\begin{align*}
	\beta_\delta := 
	\inf_{u_\delta\in\U_\delta} \sup_{v_\delta\in\V_\delta} 
		\frac{b(u_\delta,v_\delta)}{\| u_\delta\|_\U\, \| v_\delta\|_\V} = 1.
\end{align*}
\end{prop}
\begin{proof}
	Let $0\ne u_\delta\in\U_\delta\subset L_2(I;H)$. Then, since $\U_\delta=B^*(\V_\delta)$ there exists a unique $z_\delta\in\V_\delta$ such that $B^*z_\delta=u_\delta$. Hence
	\begin{align*}
		\sup_{v_\delta\in\V_\delta}  \frac{b(u_\delta,v_\delta)}{\| u_\delta\|_\U\, \| v_\delta\|_\V}
		&\ge \frac{b(u_\delta,z_\delta)}{\| u_\delta\|_\U\, \| z_\delta\|_\V}
		= \frac{(u_\delta,B^*z_\delta)_\cH}{\| u_\delta\|_\U\, \| z_\delta\|_\V}
		= \frac{(u_\delta,u_\delta)_\cH}{\| u_\delta\|_\cH\, \| B^*z_\delta\|_\cH} \\
		&=\frac{ \|u_\delta\|^2_\cH}{\| u_\delta\|_\cH\, \| u_\delta\|_\cH} = 1.
	\end{align*}
	On the other hand, by the Cauchy-Schwarz inequality, we have
	\begin{align*}
		\sup_{v_\delta\in\V_\delta}  \frac{b(u_\delta,v_\delta)}{\| u_\delta\|_\U\, \| v_\delta\|_\V}
		&= \sup_{v_\delta\in\V_\delta}  \frac{(u_\delta, B^*v_\delta)_\cH}{\| u_\delta\|_\U\, \| v_\delta\|_\V}
		\le \sup_{v_\delta\in\V_\delta}  \frac{ \| u_\delta\|_\cH\, \|B^*v_\delta\|_\cH}{\| u_\delta\|_\cH\, \| B^*v_\delta\|_\cH} = 1,
	\end{align*}
	which proves the claim.
\end{proof}

%----------------------------------------------------------------------------------------------------------
\subsection{Optimal very weak discretization of ordinary differential equations}
%----------------------------------------------------------------------------------------------------------
For the understanding of our subsequent numerical investigations, it is worth considering the univariate case, i.e., ordinary differential equations (ODEs) of the form
\begin{align}\label{Eq:ODE}
	-u''(x) &= f(x),
	\qquad x\in (0,1),
\end{align}
with either boundary or second order initial conditions, namely
\begin{subequations}\label{BC}
	\begin{align}
		u(0) &= u(1) = 0 \quad\text{or } \label{BC:1}\\
		u(0) &=0, \quad u'(0)=0. \label{BC:2}
	\end{align}
\end{subequations}
Using the above framework, we obtain $b(u,v):= -(u,v'')_{L_2(0,1)}$ and $\U:=L_2(0,1)$ in both cases. Moreover, in this univariate setting, we can identify the test space $\V$ given in \eqref{Def:V} as follows
\begin{subequations}
	\begin{align}
		\V_{\text{BVP}} &:= H^1_0(0,1) \cap H^2(0,1),  		& \text{for }\eqref{BC:1},\\
		\V_{\text{IVP}} &:= H^2_{\{T\}}(0,1)				& \text{for }\eqref{BC:2},
	\end{align}
\end{subequations}
where $H^2_{\{T\}}(0,1)$ is defined after \eqref{eq:Vcirc}. Hence, in the ODE case, we get $\V_\circ=\V$, which makes the discretization particularly straightforward. In fact, we use B-spline bases of different orders $r\ge 1$ (i.e., polynomial degree $r-1$). The boundary conditions for \eqref{BC} can be realized by multiple knots and then omitting those B-splines at the boundaries which do not satisfy the particular boundary condition, see again Figure \ref{1D:Examples}.

\begin{wrapfigure}{l}{0.45\textwidth}
\begin{center}
%-------------------------
\begin{tikzpicture}[xscale=0.55,yscale=0.55]
	\begin{loglogaxis}[xlabel={$\#$ d.o.f},
			ylabel={$\| u-u_h\|_{L_2(I)}$},
        			title={Error and condition number over d.o.f. },
        			legend style={at={(0.13,0.67)},anchor=north,legend cell align=left} %
        			]
	\addplot+ [thick,mark=x,blue] table [x=dof, y=l2error] {Data/Ex3_Order03_oA02_oT02.dat};
	\addlegendentry{$1^*/3$}
	\addplot+ [thick,mark=*,red, mark options={fill=white}] table [x=dof, y=l2error] {Data/Ex3_Order13_oA00_oT02.dat};
	\addlegendentry{$1/3$}
	\addplot+ [thick,mark=triangle*,black,mark options={fill=white}] table [x=dof, y=l2error] {Data/Ex3_Order24_oA00_oT02.dat};
	\addlegendentry{$2/4$}
	\end{loglogaxis}
	\begin{loglogaxis}[xlabel={$\#$ d.o.f},
			ylabel={$\kappa_2(B)$},
			hide x axis,
			hide y axis]
	\addplot+ [thick,dashed,mark=x,blue,mark options=solid] table [x=dof, y=cond] {Data/Ex3_Order03_oA02_oT02.dat};
	\addplot+ [thick,dashed,mark=*,red, mark options={solid,fill=white}] table [x=dof, y=cond] {Data/Ex3_Order13_oA00_oT02.dat};
	\addplot+ [thick,dashed,mark=triangle*,black,mark options={solid,fill=white}] table [x=dof, y=cond] {Data/Ex3_Order24_oA00_oT02.dat};
	%-----------------
	\end{loglogaxis}
   	\pgfplotsset{every axis y label/.append style={rotate=180,yshift=0.1cm}}
   	\begin{loglogaxis}[
	        hide x axis,
	        ymin = 1,
	        ymax = 1e12,
	        ytick= {1e0,1e3,1e6,1e9,1e12},
	        yticklabels = {$10^0$,$10^3$,$10^6$,$10^9$,$10^{12}$},
         	axis y line*=right,
  		ylabel={$\kappa_2(B)$},
     		]
   	\end{loglogaxis}
\end{tikzpicture}
%-------------------------
\caption{\label{Fig:1DBVP}
Initial value problem \eqref{BC:2}. B-spline discretization of order $r_{\text{ansatz}}/r_{\text{test}}$, where ${}^*$ means that $X_\Dt=B^*(Y_\Dt)$.}
\end{center}
\end{wrapfigure}
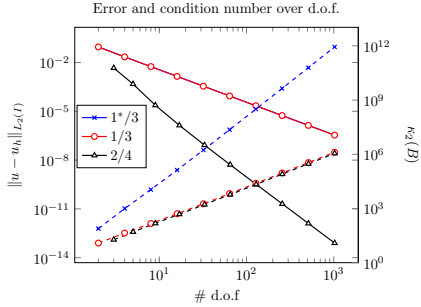

We did experiments for a whole variety of problems admitting solutions of different smoothness. The effect was negligible as we can also deduce from the left graph in Figure \ref{Fig:1DBVP}, where we depict the error and the condition number for the initial value problem \eqref{BC:2}. The three shown discretizations include the inf-sup-optimal one and test order $3$ (denoted by $1^*/3$), \enquote{standard} constant/quadratic ($1/3$), and linear/cubic ($2/4$) splines. We see the expected higher order of convergence for $2/4$ for both examples. Concerning the condition numbers, we obtain the expected $h^{-2}$ for the lower order and $h^{-4}$ for the higher order discretizations.

It is worth mentioning that we got $\beta_{\Delta t}\equiv 1$ in \emph{all} cases. This means in particular that the ansatz spaces generated by the inf-sup-optimal setting $1^*/3$ are identical with those for the $1/3$ case. After observing this numerically, we have also proven this observation. However, we stress the fact that this is a pure univariate fact, i.e., for the ODE. It is no longer true in the PDE case as we shall also see below.

%===============================================
\section{Derivation and Properties of the Algebraic Linear System}
\label{Sec:4}
%===============================================

%----------------------------------------------------------------------------------
\subsection{The linear system}
\label{SubSec:LinSystem}
%----------------------------------------------------------------------------------
To derive the stiffness matrix, we first use arbitrary spaces induced by $\{\psi_\mu := \sigma^\ell\otimes\xi_j :\, \mu=1,...,\cN_\delta\}$ for the trial and $\{  \varphi_\nu = \varrho^k\otimes\phi_i:\, \nu=1,...,\cN_\delta\}$ for the test space. Using $[\bbB_\delta]_{\mu,\nu}=[\bbB_\delta]_{(\ell,j),(k,i)}$ we get
\begin{align}
		 [\bbB_\delta]_{(\ell,j),(k,i)} 
		&= b(\psi_\mu,\varphi_\nu) 
		= (\psi_\mu,B^*\varphi_\nu)_\cH  %\nonumber \\
		 = ( \sigma^\ell\otimes\xi_j,
			\ddot{\varrho}^k \otimes \phi_i + \varrho^k\otimes  A\phi_i)_\cH \nonumber\\
		&= ( \sigma^\ell, \ddot{\varrho}^k)_{L_2(I)}\, (\xi_j, \phi_i)_{L_2(\Omega)}
		+ ( \sigma^\ell, {\varrho}^k)_{L_2(I)}\, (\xi_j, A\phi_i)_{L_2(\Omega)},
		\label{Eq:Bdelta-general}
\end{align}
so that $\bbB_\delta= \tilde\bN_\Dt \otimes \tilde\bM_h + \tilde\bM_\Dt\otimes\tilde\bN_h$, where 
$[\tilde\bM_\Dt]_{\ell,k} := ({\sigma}^\ell, {\varrho}^k)_{L_2(I)}$, 
$[\tilde\bM_h]_{j,i} := (\xi_j, \phi_i)_{L_2(\Omega)}$, 
$[\tilde\bN_\Dt]_{\ell,k} := ({\sigma}^\ell, \ddot{\varrho}^k)_{L_2(I)}$ 
and $[\tilde\bN_h]_{j,i} := (\xi_j, A\phi_i)_{L_2(\Omega)}$. 
In the specific case $\psi_\mu=B^*(\varphi_\mu)$, we get the representation
\begin{align}
	 [\bbB_\delta]_{(\ell,j),(k,i)} 
	&= b(\psi_\mu,\varphi_\nu) 
		= (\psi_\mu,B^*\varphi_\nu)_\cH  
		= (B^*\varphi_\mu,B^*\varphi_\nu)_\cH 
					\nonumber\\
		&= (\ddot{\varrho}^\ell \otimes \phi_j + \varrho^\ell\otimes  A\phi_j,
			\ddot{\varrho}^k \otimes \phi_i + \varrho^k\otimes  A\phi_i)_\cH
			\nonumber\\
		&= (\ddot{\varrho}^\ell, \ddot{\varrho}^k)_{L_2(I)}
			\, (\phi_j, \phi_i)_{L_2(\Omega)}
			+ ({\varrho}^\ell, {\varrho}^k)_{L_2(I)}
			\, (A\phi_j, A\phi_i)_{L_2(\Omega)}
			\nonumber\\
		&\qquad  + (\ddot{\varrho}^\ell, \varrho^k)_{L_2(I)}
			\, (\phi_j, A\phi_i)_{L_2(\Omega)}
			+ ({\varrho}^\ell, \ddot{\varrho}^k)_{L_2(I)}
			\, (A\phi_j, \phi_i)_{L_2(\Omega)}
					\label{Eq:Bdelta-optimal}
\end{align}
so that $\bbB_\delta= \bQ_\Dt \otimes \bM_h + \bN_\Dt\otimes \bN_h^\top +\bN_\Dt^\top\otimes \bN_h + \bM_\Dt\otimes\bQ_h$, where
\begin{align*}
	[\bQ_\Dt]_{\ell,k} &:= (\ddot{\varrho}^\ell, \ddot{\varrho}^k)_{L_2(I)},
	&\kern-9pt
	[\bM_\Dt]_{\ell,k} &:= ({\varrho}^\ell, {\varrho}^k)_{L_2(I)},
	&\kern-3pt
	[\bN_\Dt]_{\ell,k} &:= (\ddot{\varrho}^\ell, {\varrho}^k)_{L_2(I)},\\
	%--
	[\bQ_h]_{j,i} &:= (A\phi_j, A\phi_i)_{L_2(\Omega)}, &
	[\bM_h]_{j,i} &:= (\phi_j, \phi_i)_{L_2(\Omega)}, &
	[\bN_h]_{j,i} &:= (A\phi_j, \phi_i)_{L_2(\Omega)}.
\end{align*}
We stress that $\bbB_\delta$ is symmetric and positive definite for $A=-\Delta$. 
Finally, let us now detail the right-hand side. Recall 
from \eqref{eq:veryweak}, that $g(v)= (f,v)_\cH + \langle u_1, v(0)\rangle - (u_0, \dot{v}(0))_H$. Hence,
\begin{align*}
	[\bg_\delta]_\nu 
		&= [\bg_\delta]_{(k,i)}
		= (f,\varphi_\nu)_\cH + \langle u_1, \varphi_\nu(0)\rangle_{V'\times V} - (u_0, \dot{\varphi_\nu}(0))_H \\
		&= (f, \varrho^k\otimes \phi_i)_\cH + \langle u_1, \varphi_\nu(0)\rangle_{V'\times V} - (u_0, \dot{\varphi_\nu}(0))_H \\
		&= \int_0^T \int_\Omega f(t,x)\, \varrho^k(t)\, \phi_i(x)\, dx\, dt
			+ \int_\Omega [ u_1(x)\, \varrho^k(0) - u_0(x)\,  \dot\varrho^k(0)]   \phi_i(x)\, dx.
\end{align*}
Using appropriate quadrature formulae results in a numerical approximation, which we will again denote by $\bg_\delta$. 
Then, solving the linear system $\bbB_\delta \bu_\delta=\bg_\delta$ yields the expansion coefficients of the desired approximation $u_\delta\in\U_\delta$ as follows: Let $\bu_\delta=(u_\mu)_{\mu=1,...,\cN_\delta}$, $\mu=(k,i)$, then
\begin{align*}
	u_\delta(t,x)
	&= \sum_{\mu=1}^{\cN_\delta}  u_\mu\, \psi_\mu(t,x)
	= \sum_{k=1}^{N_t}\sum_{i=1}^{N_h} u_{k,i}\, \sigma^k(x)\, \xi_i(x), 
\end{align*}
in the general case and for the special one, i.e., $\psi_\mu=B^*(\varphi_\mu)$,
\begin{align*}
	u_\delta(t,x)
	&= \sum_{\mu=1}^{\cN_\delta}  u_\mu\, \psi_\mu(t,x)
	= \sum_{k=1}^{N_t}\sum_{i=1}^{N_h} u_{k,i}\, 
		\big(\ddot{\varrho}^k(t)\, \phi_i(x) + \varrho^k(t)\,  A\phi_i(x)\big).
\end{align*}

%----------------------------------------------------------------------------------
\subsection{Stability vs.\ conditioning}
\label{sec:stab}
%----------------------------------------------------------------------------------
The (discrete) inf-sup constant refers to the stability of the discrete system, being included in the error/residual relation
\begin{align*}
	\| u^*-u^*_\delta\|_\U \le \frac{1}{\beta} \sup_{v\in\V} \frac{g(v) - b(u^*_\delta,v)}{\| v\|_\V} = \frac{1}{\beta} \| r_\delta\|_{\V'},
\end{align*}
where the residual $r_\delta\in\V'$ is defined as usual by $r_\delta(v):=g(v) - b(u^*_\delta,v)$, $v\in\V$. The inf-sup constant is  the minimal generalized eigenvalue of a generalized eigenvalue problem and its continuous analogue, respectively. This has no effect on the condition number $\kappa(\bbB_\delta)$, which instead
governs the accuracy of direct solvers and convergence of iterative methods in the symmetric case.

\subsubsection*{Conditioning of the matrices}
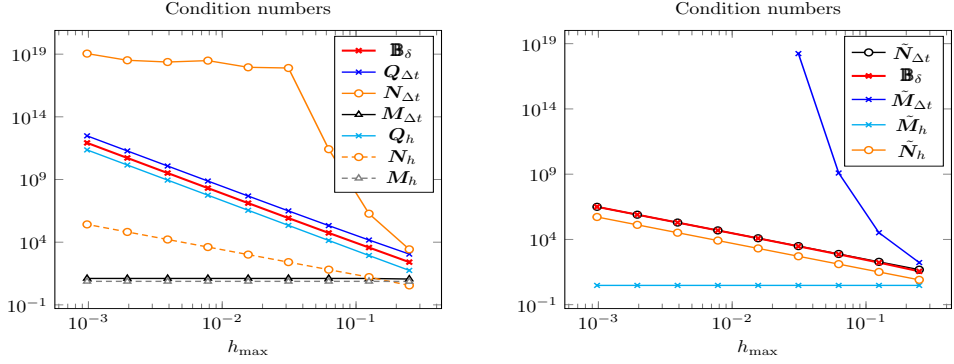
\begin{figure}[!htb]
	\begin{center}
\begin{tikzpicture}[xscale=0.75,yscale=0.65]
	\begin{loglogaxis}[xlabel={$h_{\max}$},
        			title={Condition numbers}]
	%---
	\addplot+ [very thick,mark=x,red,mark options={solid,fill=white}] table [x=hmax, y=condB] {Data/resultsWaveCondOpt_03_03.txt};
	\addlegendentry{$\bbB_\delta$}
	%---
	\addplot+ [thick,mark=x,blue] table [x=hmax, y=condQtime] {Data/resultsWaveCondOpt_03_03.txt};
	\addlegendentry{$\bQ_\Dt$}
	\addplot+ [thick,mark=*,orange, mark options={fill=white}] table [x=hmax, y=condDtime] {Data/resultsWaveCondOpt_03_03.txt};
	\addlegendentry{$\bN_\Dt$}
	\addplot+ [thick,mark=triangle*,black,mark options={fill=white}] table [x=hmax, y=condMtime] {Data/resultsWaveCondOpt_03_03.txt};
	\addlegendentry{$\bM_\Dt$}
	%---
	\addplot+ [thick,mark=x,cyan] table [x=hmax, y=condQspace] {Data/resultsWaveCondOpt_03_03.txt};
	\addlegendentry{$\bQ_h$}
	\addplot+ [thick,mark=*,orange, mark options={solid,fill=white}] table [x=hmax, y=condAspace] {Data/resultsWaveCondOpt_03_03.txt};
	\addlegendentry{$\bN_h$}
	\addplot+ [thick,mark=triangle*,gray,mark options={solid,fill=white}] table [x=hmax, y=condMspace] {Data/resultsWaveCondOpt_03_03.txt};
	\addlegendentry{$\bM_h$}
	%---
	\end{loglogaxis}
\end{tikzpicture}
\hfill
\begin{tikzpicture}[xscale=0.75,yscale=0.65]
	\begin{loglogaxis}[xlabel={$h_{\max}$},
        			title={Condition numbers}]
	\addplot+ [thick,mark=*,black, mark options={fill=white}] table [x=hmax, y=condNtime] {Data/resultsWaveCond_13_13.txt};
	\addlegendentry{$\tilde\bN_\Dt$}
	\addplot+ [very thick,mark=x,red,mark options={solid,fill=white}] table [x=hmax, y=condB] {Data/resultsWaveCond_13_13.txt};
	\addlegendentry{$\bbB_\delta$}
	%---
	\addplot+ [thick,mark=x,blue] table [x=hmax, y=condMtime] {Data/resultsWaveCond_13_13.txt};
	\addlegendentry{$\tilde\bM_\Dt$}
	%---
	\addplot+ [thick,mark=x,cyan] table [x=hmax, y=condMspace] {Data/resultsWaveCond_13_13.txt};
	\addlegendentry{$\tilde\bM_h$}
	\addplot+ [thick,mark=*,orange, mark options={solid,fill=white}] table [x=hmax, y=condNspace] {Data/resultsWaveCond_13_13.txt};
	\addlegendentry{$\tilde\bN_h$}
	\end{loglogaxis}
\end{tikzpicture}
	\caption{\label{Fig:CondBdelta}Condition numbers of involved matrices, for the general case \eqref{Eq:Bdelta-general} (left) and the inf-sup-optimal case \eqref{Eq:Bdelta-optimal} (right). }
	\end{center}
\end{figure}

We report on the condition numbers of the matrices involved in \eqref{Eq:Bdelta-general} and \eqref{Eq:Bdelta-optimal}. In Figure \ref{Fig:CondBdelta}, we see the asymptotic behavior of the different matrices. Most matrices show a \enquote{normal} scaling in the order given by the order of the differential operator. However, there are two components, namely $\tilde{\textbf{M}}_{\Delta t}$ and $\textbf{N}_{\Delta t}$, which show a very poor scaling as the mesh size tends to zero (here indicated by $h_{\text{max}}$ but used for both $\Delta t$ and $h$).  As a result, the stiffness matrix shows an asymptotic behavior calling for structure-aware preconditioning.

\subsubsection*{Preconditioning}
Let  
$\bbM_\delta := \bM_{\Delta t} \otimes \bM_h$ and 
$\bbK_\delta := \bN_{\Delta t}\otimes \bM_h + \bM_{\Delta t}\otimes\bN_h$. 
Then 
\begin{align*}
	\bbK_\delta^\top  \bbM_\delta^{-1} \bbK_\delta
	&=  (\bN_{\Delta t}^\top \bM_{\Delta t}^{-1} \bN_{\Delta t}) \otimes \bM_h 
		+ \bN_\Dt\otimes \bN_h^\top 
		+\bN_\Dt^\top\otimes \bN_h \\
		&\quad + \bM_\Dt\otimes (\bN_{h}^\top \bM_{h}^{-1} \bN_{h}),
\end{align*}
so that $\bbK_\delta^\top \bbM_\delta^{-1} \bbK_\delta = \bbB_\delta$ if and only if $\bQ_\Dt = \bN_{\Delta t}^\top \bM_{\Delta t}^{-1} \bN_{\Delta t}$ and $\bQ_h = \bN_{h}^\top \bM_{h}^{-1} \bN_{h}$.

Even if we cannot hope that those relations hold exactly in general, we are going to describe situations in which at least spectral equivalence holds. To this end, we will closely follow \cite{BabushkaOsbornPitkaranta,SilvesterMihaijlovic} in a slightly generalized setting. {We} recall the biharmonic-type problem \eqref{eq:Biharm} along with its equivalent mixed form \eqref{eq:SPPBiharm}. Let us abbreviate $Z:=D(A)$ and let $Z_h :=\Span\{ \phi_1,...,\phi_{{N_h}}\}\subset Z$ be some discretization as in \eqref{eq:Zh}. Moreover, let  $H_h :=\Span\{ \xi_1,...,\xi_{{n_h}}\}\subset H$ be some finite-dimensional approximation space for the auxiliary variable. Then, setting
\begin{align*}
	\bM_h := \big[ (\xi_i, \xi_j)_H \big]_{i,j=1,...,n_h}
	&\quad
	\bA_h := \big[ (A\phi_k, \xi_j)_H \big]_{k=1,...,N_h, j=1,...,n_h}, 
\end{align*}
the discrete form of \eqref{eq:SPPBiharm} aims to determine $\bu_h\in\R^{n_h}$ and $\bz_h\in\R^{N_h}$ such that
\begin{align}\label{eq:SPPBiharmDisc}
	\begin{pmatrix}
		\bM_h & \bA_h \\ \bA_h^\top & \bNull 
	\end{pmatrix}
	\begin{pmatrix} \bu_h \\ \bz_h \end{pmatrix}
	= 
	\begin{pmatrix} \bNull \\ -\bg_h \end{pmatrix},
\end{align}
where $\bg_h = [ \langle g, \phi_k\rangle_{Z'\times Z} ]_{k=1,...,N_h}$. Note, that $\bM_h$ is symmetric and positive definite. 
The corresponding discrete operators are defined as follows
\begin{align*}
	A_h: Z_h \to H_h:\quad& (A_h z_h, u_h)_H := (A z_h, u_h)_H,\,\, u_h\in H_h, z_h\in Z_h, \\
	M_h: H_h \to H_h:\quad& (M_h u_h, v_h)_H := (u_h, v_h)_H,\,\, u_h, v_h\in M_h.
\end{align*}
The stiffness matrix for the biharmonic-type problem reads as follows: $\bQ_h:= [ (A\phi_k, A\phi_\ell)_H ]_{k,\ell =1,...,N_h}$. 
Finally, we define discrete norms on $Z_h$ by $\| z_h\|_{Z_h}^2 := \bz_h^\top \bQ_h \bz_h$ for $z_h = \sum_{k=1}^{N_h} (\bz_h)_k\,\phi_k\in Z_h$, $\bz_h\in\R^{N_h}$ and $\| u_h\|_{M_h}^2 = \bu_h^\top \bM_h \bu_h$ for $u_h = \sum_{i=1}^{n_h} (\bu_h)_i\,\xi_i\in M_h$, $\bu_h\in\R^{n_h}$.

\begin{prop}\label{Prop:SpectrEquiv}
	Let $A_h$ be \emph{bounded}, i.e., there exists a constant $0<\Gamma<\infty$ such that $(A_h z_h, u_h)_H\le\Gamma \|  z_h\|_{Z_h}\, \| u_h\|_{M_h}$ for all $u_h\in H_h$ and $z_h\in Z_h$, and \emph{uniformly inf-sup stable}, i.e.,
	\begin{align}\label{eq:infsup-a}
		\inf_{z_h\in Z_h} \sup_{u_h\in H_h} \frac{(A_h z_h, u_h)_H}{\| z_h\|_{Z_h}\, \|u_h\|_{M_h}}\ge \gamma >0.
	\end{align}
	Then, $\bQ_h$ and $\bA_h \bM_h^{-1} \bA_h^\top$ are spectrally equivalent, i.e., 
	\begin{align*}
		\gamma^2\,  \bz_h^\top \bQ_h \bz_h
		\le \bz_h^\top \bA_h \bM_h^{-1} \bA_h^\top \bz_h
		\le \Gamma^2\,  \bz_h^\top \bQ_h \bz_h
		\quad\text{ for all } \bz_h\in\R^{N_h}.
	\end{align*}
\end{prop}
\begin{proof}
	The  proof follows the lines in \cite[(1.9-1.12)]{SilvesterMihaijlovic}. 
	Let $Z_h\ni z_h = \sum_{k=1}^{N_h} (\bz_h)_k\,\phi_k$ and $M_h\ni u_h = \sum_{i=1}^{n_h} (\bu_h)_i\,\xi_i$. Then, by \eqref{eq:infsup-a}
\begin{align*}
	\gamma\,  (\bz_h^\top \bQ_h \bz_h)^{1/2}
	&= \gamma \| z_h\|_{Z_h}
	\le \sup_{u_h\in H_h} \frac{(A_h z_h, u_h)_H}{\|u_h\|_{M_h}} 
	= \max_{\bu_h\in\R^{n_h}} \frac{\bz_h^\top \bA_h \bu_h}{(\bu_h^\top \bM_h\bu_h)^{1/2}}\\
	&= \max_{\bv_h= \bM_h^{1/2}\bu_h\in\R^{n_h}} \frac{\bz_h^\top \bA_h \bM_h^{-1/2} \bv_h}{(\bv_h^\top \bv_h)^{1/2}}
	= (\bz_h^\top \bA_h \bM_h^{-1} \bA_h^\top \bz_h)^{1/2},
\end{align*}
since it is easily seen that the maximum is attained for $\bv_h=\bM_h^{-1/2} \bA_h^\top \bz_h$, which proves the first inequality. Using the boundedness of $A_h$ yields 
$	\Gamma\,  (\bz_h^\top \bQ_h \bz_h)^{1/2}
	= \Gamma \| z_h\|_{Z_h}
	\ge \sup_{u_h\in H_h} \frac{(A_h z_h, u_h)_H}{\|u_h\|_{M_h}}$, 
so that the second inequality follows the above lines.
\end{proof}

\begin{rem}\label{Rem:SpaceSpectralEquivalence}
	In \cite[\S4]{BabushkaOsbornPitkaranta} the above assumptions have been shown within the so-called \emph{Ciarlet-Raviart method}, where $A=-\Delta$ with homogeneous Dirichlet boundary conditions on a bounded convex polygon $\Omega\subset\R^2$. Then, $D(A)=H^2(\Omega)\cap H^1_0(\Omega)$ and $H=L_2(\Omega)$ -- exactly our setting for the wave equation.
	
	Let $\{ \cT_h\}_{0<h<1}$ be a family of shape regular and quasi uniform triangulations of $\Omega$ consisting of triangles of diameter less or equal to $h$. The next piece consists of mesh dependent norms and spaces defined as $H^2_h:=\{ u\in H^1(\Omega):\, u|_T\in H^2(T), T\in\cT_h\}$, $\Gamma_h:=\bigcup_{T\in\cT_h} \partial T$ and  $\| u\|_{2,h}^2 := \sum_{T\in\cT_h} \| u\|_{2,T}^2 + h^{-1} \int_{\Gamma_h} \left| J \frac{\partial u}{\partial\nu}\right|^2\, ds$,
	where
	\begin{align*}
		\left. J \frac{\partial u}{\partial\nu}\right|_{T'}
		&:=
		\begin{cases}
			\frac{\partial u}{\partial\nu^1}+\frac{\partial u}{\partial\nu^2},
				& \text{if }\, T'=\partial T^1\cap \partial T^2 \text{ is an interior edge of } \cT_h, \\[3pt]
			\frac{\partial u}{\partial\nu},
				& \text{if } T' \text{is a boundary edge of } \cT_h,
		\end{cases}
	\end{align*}
	and $\nu^j$ denotes the unit outward normal of $T^j$. Next, let
	\begin{align*}
		\| u\|_{0,h}^2 &:= \| u\|_{L_2(\Omega)}^2 + h\, \int_{\Gamma_h} |u(s)|^2\, ds,
		\qquad u\in H^1(\Omega)
	\end{align*}
	and define $H_h^0$ as the completion of $H^1(\Omega)$ w.r.t.\ $\| \cdot\|_{0,h}$. Then $H_h^0\cong L_2(\Omega)\oplus L_2(\Gamma_h)$. For $S_h:=\{ v\in C^0(\bar\Omega):\, v|_T\in\cP_k, T\in\cT_h\}$, $k\ge 1$ and $\cP_k$ denoting the space of polynomials of degree $k$ or less, we have that $S_h\subset H_h^0\cap H^2_h$. For $k=3$, we get that $S_h=Z_h$ with $Z_h$ defined in \eqref{eq:Zh}.
		The discrete operator $A_h$ is induced by the bilinear form $a_h: H_h^0\times H_h^2\cap H^1_0(\Omega)\to\R$ defined by
		\begin{align*}
			a_h(u_h, w_h)
			&:= \sum_{T\in\cT_h} \int_T u_h(x)\, \Delta w_h(x)\, dx - \int_{\Gamma_h} u_h\! \left( J \frac{\partial w_h}{\partial\nu}\right)\! ds
		\end{align*}
	and $M_h$ is induced by $m_h(u,v):=(u,v)_{L_2(\Omega)}$. 
	The discrete spaces arise there from a shape regular and quasi uniform triangulation of $\Omega$ as well as mesh-dependent inner products and norms. 
	The boundedness of $A_h$ is immediate. The inf-sup-stability \eqref{eq:infsup-a} was proven in \cite[Thm.\ 3]{BabushkaOsbornPitkaranta}.
	
	Noting that $a_h(u_h, w_h)=-(\nabla u_h, \nabla w_h)_{L_2(\Omega)}$ for $u\in H^1(\Omega)$ and $w_h\in H_h^2$, we obtain that $\bQ_h$ and $\bN_{h}^\top \bM_{h}^{-1} \bN_{h}$ defined in \S\ref{SubSec:LinSystem} are in fact spectrally equivalent.
\end{rem}

We observed the spectral equivalence for the spatial matrices also in our numerical experiments. However, we saw that this is not true for the temporal matrices in the sense that $\bQ_\Dt$ and $\bN_{\Delta t}^\top \bM_{\Delta t}^{-1} \bN_{\Delta t}$ are not spectrally equivalent.

%===============================================
\section{Solution of the Algebraic Linear System}
\label{Sec:5}
%===============================================
To derive preconditioning strategies and the new projection method,
we rewrite the linear system $\bbB_\delta \bu_\delta = \bg_\delta$
as a linear matrix equation, so as to exploit the {structure of the Kronecker problem}.  
Let $\bx = {\rm vec}(\bX)$ be the operator stacking the columns of $\bX$ one after the
other, then it holds that $(\bB \otimes \bA) \bx = {\rm vec}(\bA \bX \bB^\top)$ for given matrices
$\bA, \bX$, and $\bB$ of conforming dimensions. Hence, the vector system is written as
\begin{equation}\label{eq:matrix_eq}
{\mathcal A}(\bU) = \bG, \,\mbox{with}\,
{\mathcal A}(\bU)=\bM_h \bU \bQ_\Dt^\top + \bN_h^\top \bU \bN_\Dt^\top + \bN_h \bU \bN_\Dt + \bQ_h \bU \bM_\Dt,
\end{equation}
where $\bg={\rm vec}(\bG)$ and the symmetry of some of the matrices has been exploited.

In the following we describe two distinct approaches: 
First, we recall the matrix-oriented conjugate gradient method,
preconditioned by two different operator-aware  strategies. Then we discuss a procedure
that directly deals with (\ref{eq:matrix_eq}).

%----------------------------------------------------------------------------------
\subsection{Preconditioned conjugate gradients}
%----------------------------------------------------------------------------------
Since $\bbB_\delta$ is symmetric and positive definite, the preconditioned conjugate gradient (PCG) method can be  applied directly to \eqref{eq:matrix_eq}, yielding  a matrix-oriented implementation of PCG, see Algorithm \ref{Alg:MatrixPCG}. Here tr$(\bX)$ denotes the trace of the square matrix $\bX$. In exact precision arithmetic, this formulation,  gives the same iterates as the standard vector form, while exploiting matrix-matrix computations\footnote{The matrix-oriented version of PCG is also used to exploit low rank representations of the iterates, in case the starting residual is low rank and the final solution can be well approximated by a low rank matrix; see, e.g., \cite{Kressner.Tobler.11}. We will not exploit this setting here.}.

\renewcommand{\algorithmicrequire}{\textbf{Input:}}
\begin{algorithm}[t]
\caption{Matrix-oriented PCG}\label{Alg:MatrixPCG}
\begin{algorithmic}[1]
\Require $\bU_0$
\State set $\bR_0= \bG - {\mathcal A}(\bU_0)$, $\bZ_0 = {\mathcal P}^{-1}(\bR_0)$,  $\bP_0=\bZ_0$, $\gamma_0={\rm tr}(\bR_0^\top \bZ_0)$% $\beta=0$, $\bP_{-1}=0$
\For{$k=0,1,...$}
	\State $\delta = {\rm tr}(\bP_k^\top {\mathcal A}(\bP_k))$, $\alpha = \gamma_k/\delta$
	\State $\bX_{k+1} = \bX_k + \alpha \bP_k$
	\State $\bR_{k+1} = \bG - {\mathcal A}(\bX_{k+1})$
	\State $\bZ_{k+1} = {\mathcal P}^{-1}(\bR_{k+1})$
	\State $\gamma_{k+1} = {\rm tr}(\bR_{k+1}^\top \bZ_{k+1})$, $\beta = \gamma_{k+1}/\gamma_{k}$
	\State $\bP_{k+1} = \bZ_{k+1} + \beta \bP_k$
	
\EndFor
\end{algorithmic}
\end{algorithm}

\subsubsection{Sylvester operator preconditioning.}\label{PCG_1} A natural preconditioning strategy consists
of taking the leading part of the coefficient matrix, in terms of order of the differential operators.
Hence, setting ${\mathbb P} =  
\bQ_\Dt \otimes \bM_h + \bM_\Dt \otimes \bQ_h$, we have (see also \cite{HPSU20})
\begin{align*}
	\bz_{k+1} = {\mathbb P}^{-1} \br_{k+1} \quad
\Leftrightarrow \quad
	\bZ_{k+1} = {\mathcal P}^{-1}(\bR_{k+1}),
\end{align*}
with $\br_{k+1}={\rm vec}(\bR_{k+1})$ and $\bz_{k+1} = {\rm vec}(\bZ_{k+1})$. Applying ${\mathcal P}^{-1}$ corresponds to solving the generalized Sylvester equation $\bM_h \bZ \bQ_\Dt^\top +  \bQ_h \bZ \bM_\Dt = \bR_{k+1}$.  For small size problems in space, this can be carried out by means of  the Bartels-Stewart method \cite{Bartels.Stewart.72}, which entails the computation of two Schur decompositions, performed before the PCG iteration is started. 
For fine discretizations in space, iterative procedures need to be used. For these purposes, we use a Galerkin approach based on the rational Krylov subspace \cite{Druskin.Simoncini.11}, only performed on the spatial matrices; see \cite{Simoncini.survey.16} for a general discussion. A key issue is that this class of iterative methods requires the right-hand side to be low rank; we deliberately set the rank to be at most four. Hence, the Sylvester solver is applied after a rank truncation of $\bR_{k+1}$, which thus becomes part of the
preconditioning application.

\subsubsection{$\bbK_\delta^\top  \bbM_\delta^{-1} \bbK_\delta$-preconditioning.}\label{PCG_2}
To derive a preconditioner that takes full account of the coefficient matrix we employ the operator  $\bbK_\delta^\top  \bbM_\delta^{-1} \bbK_\delta$ in \S\ref{sec:stab}. Thanks to the spectral equivalence in Proposition~\ref{Prop:SpectrEquiv}, PCG applied to the resulting preconditioned operator appears to be optimal, in the sense that the number of iterations to reach the required accuracy is independent of the spatial mesh size; see Table~\ref{table:radial_example_1CG}.

In vector form this preconditioner is applied as $\bz_{k+1} = \left ( \bbK_\delta^\top  \bbM_\delta^{-1} \bbK_\delta \right )^{-1} \br_{k+1}$. However, this operation can be performed without explicitly using the Kronecker form of the involved matrices, with significant computational and memory savings.
We observe that 
$$
\bbK_\delta 
	\kern-1pt=\kern-1pt \bN_{\Delta t}\otimes \bM_h + \bM_{\Delta t}\otimes\bN_h 
	\kern-1pt=\kern-1pt ( \bN_{\Delta t} \bM_{\Delta t}^{-1} \otimes {\bI} + {\bI} \otimes \bN_h \bM_h^{-1} ) (\bM_{\Delta t} \otimes \bM_h) 
	\kern-1pt=: \widehat \bbK_\delta  \bbM_\delta.
$$
Moreover, due to the transposition properties of the Kronecker product,
$\bbK_\delta^\top = \widehat \bbK_\delta^\top \bbM_\delta $.
Hence, 
$\bbK_\delta^\top  \bbM_\delta^{-1} \bbK_\delta = 
\widehat \bbK_\delta^\top \widehat \bbK_\delta  \bbM_\delta$. 
Therefore, 
$$
\bZ_{k+1} = {\mathcal P}^{-1}(\bR_{k+1}) \,\, \Leftrightarrow \,\,
\bz_{k+1} = \bbM_\delta^{-1} \widehat \bbK_\delta^{-1} (\widehat \bbK_\delta^\top)^{-1} \br_{k+1},
$$
We next observe that the equation $(\widehat \bbK_\delta^\top) \bw = \br_{k+1}$ can be
written as the following Sylvester matrix equation
\begin{equation}\label{eq:Sylv1}
 \bW \bM_{\Delta t}^{-1} \bN_{\Delta t} +   \bM_h^{-1} \bN_h \bW  = \bR_{k+1}
\end{equation}
and analogously for $(\widehat \bbK_\delta^\top) \widehat\bw = \bw$, that is
\begin{equation}\label{eq:Sylv2}
 \widehat \bW \bM_{\Delta t}^{-1} \bN_{\Delta t} +   \bM_h^{-1} \bN_h \widehat \bW  = \bW .
\end{equation}
Finally, the preconditioned matrix is obtained as $\bZ_{k+1} = \bM_h^{-1} \widehat \bW \bM_{\Delta t}^{-1}$.

Summarizing, the application of the operator preconditioner amounts to the solution
of the two Sylvester matrix equations (\ref{eq:Sylv1})-(\ref{eq:Sylv2}), and the product
$\bZ_{k+1} = \bM_h^{-1} \widehat \bW \bM_{\Delta t}^{-1}$. The overall computational cost of
this operation depends on the cost of solving the two matrix equations. For small dimensions
in space, once again a Schur-decomposition based method can be used \cite{Bartels.Stewart.72}; we recall
here that thanks to the discretization employed, we do not expect to have large dimensions in time,
as matrices of size at most ${\mathcal O}(100)$ arise. 
Also in this case, for fine discretizations in space we use an iterative method (Galerkin) based
on the rational Krylov subspace \cite{Druskin.Simoncini.11}, only performed on the spatial matrices,
with the truncation of the corresponding right-hand side, $\bR_{k+1}$ and $\bW$, respectively,
so as to have at most rank equal to four. Allowing a larger rank did not seem to improve
the effectiveness of the preconditioner. Several implementation enhancements can be developed
to make the action of the preconditioner more efficient, since most operations are repeated
at each PCG iteration with the same matrices.

%----------------------------------------------------------------------------------
\subsection{Galerkin projection}
%----------------------------------------------------------------------------------
An alternative to PCG consists of attacking the original multi-term matrix equation directly. 
Thanks to the symmetry of $\bN_h$ we rewrite the matrix equation (\ref{eq:matrix_eq}) as
\begin{equation}\label{eq:matrix_eq1}
\bM_h \bU \bQ_\Dt^\top + \bN_h^\top \bU (\bN_\Dt^\top + \bN_\Dt) + \bQ_h \bU \bM_\Dt 
 = \bG, 
\end{equation}
with $\bG$ of low rank, that is $\bG = \bG_1 \bG_2^\top$. Consider two appropriately selected vector spaces ${\mathcal V}_k$, ${\mathcal W}_k$ of dimensions much lower than $N_h, N_t$, respectively, and let $\bV_k$, $\bW_k$ be the matrices whose orthonormal columns span the two corresponding spaces. We look for a low rank approximation of $\bU$ as $\bU_k = \bV_k \bY_k \bW_k^\top$. To determine $\bY_k$ we impose an orthogonality (Galerkin) condition on the residual 
\begin{equation}\label{eq:residual_matrix}
	\bR_k := \bG_1 \bG_2^\top - \bM_h \bU_k \bQ_\Dt^\top - \bN_h^\top \bU_k (\bN_\Dt^\top + \bN_\Dt) - \bQ_h \bU_k \bM_\Dt.
\end{equation}
with respect to the generated space pair $(\bV_k, \bW_k)$. Using the matrix Euclidean inner product, this corresponds to imposing that $\bV_k^\top \bR_k \bW_k=0$. Substituting $\bR_k$ and $\bU_k$ into this matrix equation, we obtain the following {\em reduced} matrix equation, of the same type as (\ref{eq:matrix_eq1}) but of much smaller size,
\begin{align*}
	(\bV_k^\top \bM_h \bV_k) \bY_k (\bQ_\Dt^\top \bW_k) 
		&+ (\bV_k^\top\bN_h^\top \bV_k) \bY_k (\bW_k^\top (\bN_\Dt^\top + \bN_\Dt) \bW_k) \\
		&+  (\bV_k^\top\bQ_h \bV_k)  \bY_k (\bW_k^\top \bM_\Dt \bW_k) 
		= (\bV_k^\top\bG_1) (\bG_2^\top \bW_k) .
\end{align*}
The small dimensional matrix $\bY_k$ is thus obtained by solving the Kronecker form of this  equation\footnote{To this end, Algorithm~\ref{Alg:MatrixPCG} with a preconditioning strategy similar to the ones described in Sections~\ref{PCG_1}--\ref{PCG_2} can be employed as well.}. 
The described Galerkin reduction strategy has been thorough exploited and analyzed for Sylvester equations, and more recently successfully applied to {multi-term} equations, see, e.g., \cite{Powell2017}. The key problem-dependent ingredient is the choice of the spaces ${\mathcal V}_k$, ${\mathcal W}_k$, so  that they well represent spectral information of the ``left-hand'' and ``right-hand'' matrices in (\ref{eq:matrix_eq1}).
A well established choice is (a combination of) rational Krylov subspaces \cite{Simoncini.survey.16}. More precisely, for the spatial approximation we generate the growing space range($\bV_k$) as
$$
	\widehat \bV_{k+1} 
	= [\bV_k, (\bQ_h + s_k \bM_h)^{-1} \bv_k, (\bN_h + \sqrt{s_k} \bM_h)^{-1} \bv_k], 
	\quad \bV_1 = \bG_1,
$$
where $\bv_k$ is the $k$th column of $\bV_k$, so that $\bV_{k+1}$ is obtained by orthogonalizing the new columns inserted in $\widehat \bV_{k+1}$. The matrix $\widehat \bV_{k+1}$ grows at most by two vectors at the time. For each $k$, the parameter $s_k$ can be chosen either a-priori or dynamically, with the same sign as the spectrum of $\bQ_h$ ($\bN_h$). Here $s_k$ is cheaply determined using the adaptive strategy in \cite{Druskin.Simoncini.11}. Since $\bN_h$ represents an operator of the second order, the value $\sqrt{s_k}$ resulted to be appropriate; a specific computation of the parameter associated with $\bN_h$ can also be included, at low cost. Analogously, 
$$
\widehat \bW_{k+1} = [\bW_k, (\bQ_\Dt + \ell_k \bM_\Dt)^{-1} \bw_k, ((\bN_\Dt+\bN_\Dt^\top) + \sqrt{\ell_k} \bM_\Dt)^{-1} \bw_k], \quad
\bW_1 = \bG_2,
$$
where $\bw_k$ is the $k$th column of $\bW_k$, and $\bW_{k+1}$ is obtained by orthogonalizing the
new columns inserted in $\widehat \bW_{k+1}$. The choice of $\ell_k>0$ is made as for $s_k$.

\begin{rem}
This approach yields the vector approximation  $\bu_k = (\bW_k \otimes \bV_k) \by_k$, with $\by_k={\rm vec}(\bY_k)$ that is, the approximation space range($\bW_k \otimes \bV_k$) is more structured than that generated by PCG applied to ${\mathcal A}$. Experimental evidence shows that  this structure-aware space requires significantly smaller dimension to achieve similar accuracy. This is theoretically clear in the Sylvester equation case \cite{Simoncini.survey.16}, while it is an open problem for the multi-term linear equation setting.
\end{rem}

\begin{rem}\label{rem:linearsystems_ratKrylov}
For fine space discretizations, the most expensive step of the Galerkin projection is the solution of the linear systems with $(\bQ_h + s_k \bM_h)$ and $(\bN_h + \sqrt{s_k} \bM_h)$. Depending on the size and sparsity, these systems can be solved by either a sparse direct method or by an iterative procedure; see \cite{Simoncini.survey.16} and references therein.
\end{rem}

%===============================================
\section{Numerical Experiments}
\label{Sec:6}
%===============================================
We report some results of our extensive numerical {experiments} for  the wave equation \eqref{Eq:1.1} with $A=-c^2\Delta$, $H=L_2(\Omega)$, $\Omega\subset\R^d$ some open bounded domain, $c\ne0$ being the wave speed,  $V=H^1_0(\Omega)$ and $I=(0,1)$, i.e., $T=1$.  We choose the data in such a way that the respective solutions have different regularity. In order to do so, we use $\Omega=(0,1)^d$, so that we can construct explicit solutions by the d'Alembert formula as follows. We consider rotationally symmetric problems around the center $\bc=(c_i)_{i=1,...,d}$, $c_i=0.5$. Then, we consider polar coordinates in space, i.e., $r:=\| \bx-\br\|$, $\bx\in\Omega$. For $u_1(r) \equiv f(r,t) \equiv 0$, the solution reads
\begin{align*}
	u(r,t) = (r+ct)\, u_0(r+ct) + \frac{r-ct}{2r}\, u_0(r-ct) \quad \text{for}\; r > 0.
\end{align*}
We choose $c$ in such a way that homogeneous Dirichlet conditions can be prescribed. 
If $\bU_k$ denotes the current approximate solution computed at iteration $k$,  Algorithm~\ref{Alg:MatrixPCG} and the Galerkin method are stopped as soon as the backward error $\mathcal{E}_k$ is smaller than $10^{-5}$, where $\mathcal{E}_k$ is defined as
$$
	\mathcal{E}_k=
	\frac{\|\bR_k\|_F}%
	{\|\bG\|_F+\|\bU_k\|_F(\|\bM_h\|_F\|\bQ_{\Delta t}\|_F+\|\bQ_h\|_F\|\bM_{\Delta t}\|_F+2\|\bN_h\|_F\|\bN_{\Delta t}\|_F)},
$$
where $\bR_k$ is the residual matrix defined in~\eqref{eq:residual_matrix}.  For the Galerkin approach the computation of $\mathcal{E}_k$ simplifies thanks to the low-rank format of  the involved quantities (for instance, $\bR_k$ does not need to be explicitly formed to compute its norm). Moreover, the linear systems in the rational Krylov subspace basis construction are solved by the vector PCG method with a tolerance $\epsilon=10^{-8}$; see Remark~\ref{rem:linearsystems_ratKrylov}.

We compared the space-time method with the classical Crank-Nicolson time stepping scheme, in terms of approximation accuracy and CPU time. The $N_h\times N_h$ linear systems involved in the time marching scheme are solved by means of the vector PCG method with tolerance $\epsilon=10^{-6}$.

The code is ran in \texttt{Matlab} and the B-spline implementation is based on \cite{Mollet16}\footnote{Executed on the BwUniCluster 2.0 on instances with 32GB of RAM on two cores of an Intel Xeon Gold 6230.}. To explore the potential of the new very weak method on low-regularity solutions, we only concentrate on experiments with lower regularity solutions, in particular a solution which is continuous with discontinuous derivative (Case 1) and a discontinuous solution (Case 2). This is realized through the choice of $u_0$. On the other hand, for smooth solutions the time-stepping method would be expected to be more accurate, due to its second-order convergence, compared to the very weak method, as long as the latter uses piecewise constant trial functions.

We describe our results for the 3D setting, with $\Omega=(0,1)^3$. The data are summarized as follows
\begin{center}
\begin{tabular}{l|l|l|}
			& Case 1 & Case 2 \\ \hline
$u_0(r)$	& $(1-5 r) \mathds{1}_{r<0.2}$			&	$\mathds{1}_{r<0.2}$   \\
$c$		& $0.2$							&	$0.2$	\\
$u$		& $\in C( \bar{I} \times \bar\Omega)\setminus C^1( {I} \times\Omega)$ 	&	$\not\in C( \bar{I} \times \bar\Omega)$
\end{tabular}
\end{center}
We use tensor product spaces for the spatial discretization for both approaches. In the space-time setting we use B-splines in each direction for the test functions. For the time-stepping method, we use a Galerkin approach in which the trial and test functions are given by B-splines. Hence, the radial symmetry cannot be exploited by either methods, and the tensor product approach provides no limitation. All tables show the matrix dimensions {$N_t$} in the time space (``Time'') and {$N_h$} in the spatial space (``Space''). {We display results for uniform discretizations in space and time, where $N_h=N_t^3$, but stress the fact that our space-time discretization is \emph{unconditionally} stable, i.e., for any combination of $N_t$ and $N_h$.}

%----------------------------------------------------------------------------------
\subsection{Case 1: Continuous, but not continuously differentiable solution}
%----------------------------------------------------------------------------------

We start by comparing the performance of the two preconditioners in Sections~\ref{PCG_1}--\ref{PCG_2}, namely the Sylvester operator and the $\bbK_\delta^\mathsf{T}\bbM_\delta^{-1} \bbK_\delta$ operator. 
Both preconditioners are applied inexactly as described in the corresponding sections.

The $L_2$-error, the number of iterations and the {wall-clock time (using the Matlab tic-toc commands)} are displayed in Table~\ref{table:radial_example_1CG}. We obtain comparable errors but significantly smaller CPU times for  $\bbK_\delta^\mathsf{T}\bbM_\delta^{-1} \bbK_\delta$. Moreover, the number of iterations performed by using the $\bbK_\delta^\mathsf{T}\bbM_\delta^{-1} \bbK_\delta$ preconditioner are independent of the discretization level, illustrating the spectral equivalence of Proposition~\ref{Prop:SpectrEquiv}. We observe that PCG could not be used for further refinements, due to memory constraints of 32 GB RAM.
 
\renewcommand{\arraystretch}{1.3} 
\begin{table}[htb]
\small

% Load the data
\pgfplotstableread{Data/3DRadialExample1-CG-opt-exact-0-maxIt-1000-tolerance-1e-05-toleranceRes-0.01.txt}\ThreeDRadialOneCGOptFalse
%\pgfplotstableread{data/3DRadialExample1-CG-opt-exact-1-maxIt-1000-tolerance-0.0001-toleranceRes-0.1.txt}\ThreeDRadialOneCGOptTrue
\pgfplotstableread{Data/3DRadialExample1-CG-lyap-exact-0-maxIt-1000-tolerance-1e-05-toleranceRes-0.01.txt}\ThreeDRadialOneCGLyapFalse
%\pgfplotstableread{data/3DRadialExample1-CG-lyap-exact-1-maxIt-1000-tolerance-0.0001-toleranceRes-0.1.txt}\ThreeDRadialOneCGLyapTrue
\pgfplotstableread{Data/3DRadialExample1-Galerkin-0-maxIt-1000-tolerance-1e-05-toleranceRes-0.01.txt}\ThreeDRadialOneGalerkin
\pgfplotstableread{Data/3DRadialExample1-TS--tolerance-1e-06.txt}\ThreeDRadialOneTimestepping
\pgfplotstableread{Data/3DrefinementsPCG.txt}\ThreeDRadialOneCG % CG methods
\pgfplotstableread{Data/3DrefinementsPCG.txt}\ThreeDRadialOneGT % Galerkin and time stepping

% Concatinate the columns in the right order
%\pgfplotstablecreatecol[copy column from table={\ThreeDRadialOneCGOptTrue}{[index] 3}] {errorCGOPTexact} {\ThreeDRadialOne}
%\pgfplotstablecreatecol[copy column from table={\ThreeDRadialOneCGOptTrue}{[index] 2}] {timeCGOPTexact} {\ThreeDRadialOne}

\pgfplotstablecreatecol[copy column from table={\ThreeDRadialOneCGOptFalse}{[index] 3}] {errorCGOPTinexact} {\ThreeDRadialOneCG}
\pgfplotstablecreatecol[copy column from table={\ThreeDRadialOneCGOptFalse}{[index] 2}] {timeCGOPTinexact} {\ThreeDRadialOneCG}
\pgfplotstablecreatecol[copy column from table={\ThreeDRadialOneCGOptFalse}{[index] 1}] {iterCGOPTinexact} {\ThreeDRadialOneCG}

%\pgfplotstablecreatecol[copy column from table={\ThreeDRadialOneCGLyapTrue}{[index] 3}] {errorCGLYAPexact} {\ThreeDRadialOne}
%\pgfplotstablecreatecol[copy column from table={\ThreeDRadialOneCGLyapTrue}{[index] 2}] {timeCGLYAPexact} {\ThreeDRadialOne}

\pgfplotstablecreatecol[copy column from table={\ThreeDRadialOneCGLyapFalse}{[index] 3}] {errorCGLYAPinexact} {\ThreeDRadialOneCG}
\pgfplotstablecreatecol[copy column from table={\ThreeDRadialOneCGLyapFalse}{[index] 2}] {timeCGLYAPinexact} {\ThreeDRadialOneCG}
\pgfplotstablecreatecol[copy column from table={\ThreeDRadialOneCGLyapFalse}{[index] 1}] {iterCGLYAPinexact} {\ThreeDRadialOneCG}

\hspace{-1cm}

% Table with the CG Values
\pgfplotstabletypeset[
%Rename the columns
%columns/Refinements/.style={column name={Refinement}},
%columns/errorCGOPTexact/.style={column name={$L_2$ error}},
%columns/timeCGOPTexact/.style={column name={Time}},
columns/errorCGOPTinexact/.style={column name={$L_2$-error}, column type = {|r}, sci},
columns/timeCGOPTinexact/.style={column name={Wall time [s]}, column type = {|r}, sci},
columns/iterCGOPTinexact/.style={column name={Iter.}, column type = {|r}, fixed},
%columns/errorCGLYAPexact/.style={column name={$L_2$ error}},
%columns/timeCGLYAPexact/.style={column name={Time}},
columns/errorCGLYAPinexact/.style={column name={$L_2$-error}, column type = {|r}, sci},
columns/timeCGLYAPinexact/.style={column name={Wall time [s]}, column type = {|r}, sci},
columns/iterCGLYAPinexact/.style={column name={Iter.}, column type = {|r|}, fixed},
columns/Time/.style={column type = {|r}, fixed},
columns/Space/.style={column type = {|r}, fixed},
precision=2,
skip rows between index={0}{1},
% Head row		
every head row/.style={
	before row={
		\hline
		  \multicolumn{2}{|c|}{Unknowns} 
		 %& \multicolumn{2}{c|}{CG opt (exact)}
		 & \multicolumn{3}{c|}{PCG ($\bbK_\delta^\mathsf{T}\bbM_\delta^{-1} \bbK_\delta$)}
		% & \multicolumn{2}{c|}{CG lyap (exact)}
		 & \multicolumn{3}{c|}{PCG (Sylvester)} \\
	},
	after row=\hline,		
},
column type/.add={|}{},
every last row/.style={after row=\hline}	
]{\ThreeDRadialOneCG}

\caption{Case 1:  $L_2$-error, iterations and CPU time of PCG and the two proposed preconditioners.}
\label{table:radial_example_1CG}
\end{table}

\begin{table}[htb]
\small

% Load the Data
\pgfplotstableread{Data/3DRadialExample1-CG-opt-exact-0-maxIt-1000-tolerance-1e-05-toleranceRes-0.01.txt}\ThreeDRadialOneCGOptFalse
%\pgfplotstableread{Data/3DRadialExample1-CG-opt-exact-1-maxIt-1000-tolerance-0.0001-toleranceRes-0.1.txt}\ThreeDRadialOneCGOptTrue
\pgfplotstableread{Data/3DRadialExample1-CG-lyap-exact-0-maxIt-1000-tolerance-1e-05-toleranceRes-0.01.txt}\ThreeDRadialOneCGLyapFalse
%\pgfplotstableread{Data/3DRadialExample1-CG-lyap-exact-1-maxIt-1000-tolerance-0.0001-toleranceRes-0.1.txt}\ThreeDRadialOneCGLyapTrue
\pgfplotstableread{Data/3DRadialExample1-Galerkin-0-maxIt-1000-tolerance-1e-05-toleranceRes-0.01.txt}\ThreeDRadialOneGalerkin
\pgfplotstableread{Data/3DRadialExample1-TS--tolerance-1e-06.txt}\ThreeDRadialOneTimestepping
\pgfplotstableread{Data/3Drefinements.txt}\ThreeDRadialOneCG % CG methods
\pgfplotstableread{Data/3Drefinements.txt}\ThreeDRadialOneGT % Galerkin and time stepping

% Concatinate the columns in the right order

\pgfplotstablecreatecol[copy column from table={\ThreeDRadialOneGalerkin}{[index] 3}] {errorGalerkin} {\ThreeDRadialOneGT}
\pgfplotstablecreatecol[copy column from table={\ThreeDRadialOneGalerkin}{[index] 2}] {timeGalerkin} {\ThreeDRadialOneGT}
\pgfplotstablecreatecol[copy column from table={\ThreeDRadialOneGalerkin}{[index] 1}] {iterGalerkin} {\ThreeDRadialOneGT}

\pgfplotstablecreatecol[copy column from table={\ThreeDRadialOneTimestepping}{[index] 2}] {errorTimestepping} {\ThreeDRadialOneGT}
\pgfplotstablecreatecol[copy column from table={\ThreeDRadialOneTimestepping}{[index] 1}] {timeTimestepping} {\ThreeDRadialOneGT}

\hspace{-1cm}

% Table with the Galerkin and time stepping values
\pgfplotstabletypeset[
%Rename the columns
%columns/Refinements/.style={column name={Refinement}},
%columns/errorCGOPTexact/.style={column name={$L_2$ error}},
%columns/timeCGOPTexact/.style={column name={Time}},
columns/errorCGOPTinexact/.style={column name={$L_2$ error}, column type = {|r}, sci},
columns/timeCGOPTinexact/.style={column name={Wall time}, column type = {|r}, sci},
columns/iterCGOPTinexact/.style={column name={Iter.}, column type = {|r}, fixed},
%columns/errorCGLYAPexact/.style={column name={$L_2$ error}},
%columns/timeCGLYAPexact/.style={column name={Time}},
columns/errorCGLYAPinexact/.style={column name={$L_2$ error}, column type = {|r}, sci},
columns/timeCGLYAPinexact/.style={column name={Wall time [s]}, column type = {|r}, sci},
columns/iterCGLYAPinexact/.style={column name={Iter.}, column type = {|r}, fixed},
columns/errorGalerkin/.style={column name={$L_2$-error}, column type = {|r}, sci},
columns/timeGalerkin/.style={column name={Wall time [s]}, column type = {|r}, sci},
columns/iterGalerkin/.style={column name={Iter.}, column type = {|r}, fixed},
columns/errorTimestepping/.style={column name={$L_2$-error}, column type = {|r}, sci},
columns/timeTimestepping/.style={column name={Wall time},column type/.add={|r}{|}, sci},
columns/Time/.style={column type = {|r},fixed},
columns/Space/.style={column type = {|r},fixed},
precision=2,
skip rows between index={0}{1},
% Head row		
every head row/.style={
	before row={
		\hline
		  \multicolumn{2}{|c|}{Unknowns} 
		 %& \multicolumn{2}{c|}{CG opt (exact)}
		 & \multicolumn{3}{c|}{Galerkin}
		 & \multicolumn{2}{c|}{Time stepping}\\
	},
	after row=\hline,		
},
column type/.add={|}{},
every last row/.style={after row=\hline},
]{\ThreeDRadialOneGT}

\caption{Case 1:  $L_2$-error and CPU time for the Galerkin projection and the time-stepping method.}
\label{table:radial_example_1GT}
\end{table}
\normalsize

Table \ref{table:radial_example_1GT} contains the experimental results for the Galerkin projection and the time-stepping method.  Compared with Table~\ref{table:radial_example_1CG}, we clearly see that the Galerkin method outperforms both preconditioners. Furthermore, the projection method can effectively solve the problem for a further refinement level, we thus limit reporting our subsequent results to the Galerkin approach.

Let us now focus on the comparison between the Galerkin space-time method and the time-stepping approach. The {wall-clock} times of the two approaches are similar, while the $L_2$-error is greatly in favor of the space-time method. In particular, for the Galerkin method, the convergence rate is around 0.29\footnote{For piecewise constants, we expect a rate of $1/(2d)=0.16$ for $d=3$.}, whereas the time stepping method does not converge in the last step. 

%----------------------------------------------------------------------------------
\subsection{Case 2: Discontinuous}
%----------------------------------------------------------------------------------
For the case of a discontinuous solution, our results are shown in Table \ref{table:radial_example_2GT}. As in the first case, the {wall-clock} times are comparable. However, the space-time method errors are by a factor of 4 smaller than for the time marching scheme. The Galerkin method has a convergence rate of approximately 0.09 in the last step and 0.24 in the penultimate step, whereas the time-stepping method has a convergence rate of around 0.11 in the last step and 0.19 in the penultimate step. 
\begin{table}[h]
\small

% Load the data
\pgfplotstableread{Data/3DRadialExample2-CG-opt-exact-0-maxIt-1000-tolerance-1e-05-toleranceRes-0.01.txt}\ThreeDRadialTwoCGOptFalse
%\pgfplotstableread{Data/3DRadialExample2-CG-opt-exact-1-maxIt-1000-tolerance-0.0001-toleranceRes-0.1.txt}\ThreeDRadialTwoCGOptTrue
\pgfplotstableread{Data/3DRadialExample2-CG-lyap-exact-0-maxIt-1000-tolerance-1e-05-toleranceRes-0.01.txt}\ThreeDRadialTwoCGLyapFalse
%\pgfplotstableread{Data/3DRadialExample2-CG-lyap-exact-1-maxIt-1000-tolerance-0.0001-toleranceRes-0.1.txt}\ThreeDRadialTwoCGLyapTrue
\pgfplotstableread{Data/3DRadialExample2-Galerkin-0-maxIt-1000-tolerance-1e-05-toleranceRes-0.01.txt}\ThreeDRadialTwoGalerkin
\pgfplotstableread{Data/3DRadialExample2-TS--tolerance-1e-06.txt}\ThreeDRadialTwoTimestepping
\pgfplotstableread{Data/3Drefinements.txt}\ThreeDRadialTwoCG
\pgfplotstableread{Data/3Drefinements.txt}\ThreeDRadialTwoGT

% Concatinate the columns in the right order
%\pgfplotstablecreatecol[copy column from table={\ThreeDRadialTwoCGOptTrue}{[index] 3}] {errorCGOPTexact} {\ThreeDRadialTwo}
%\pgfplotstablecreatecol[copy column from table={\ThreeDRadialTwoCGOptTrue}{[index] 2}] {timeCGOPTexact} {\ThreeDRadialTwo}

\pgfplotstablecreatecol[copy column from table={\ThreeDRadialTwoCGOptFalse}{[index] 3}] {errorCGOPTinexact} {\ThreeDRadialTwoCG}
\pgfplotstablecreatecol[copy column from table={\ThreeDRadialTwoCGOptFalse}{[index] 2}] {timeCGOPTinexact} {\ThreeDRadialTwoCG}
\pgfplotstablecreatecol[copy column from table={\ThreeDRadialTwoCGOptFalse}{[index] 1}] {iterCGOPTinexact} {\ThreeDRadialTwoCG}

%\pgfplotstablecreatecol[copy column from table={\ThreeDRadialTwoCGLyapTrue}{[index] 3}] {errorCGLYAPexact} {\ThreeDRadialTwo}
%\pgfplotstablecreatecol[copy column from table={\ThreeDRadialTwoCGLyapTrue}{[index] 2}] {timeCGLYAPexact} {\ThreeDRadialTwo}

\pgfplotstablecreatecol[copy column from table={\ThreeDRadialTwoCGLyapFalse}{[index] 3}] {errorCGLYAPinexact} {\ThreeDRadialTwoCG}
\pgfplotstablecreatecol[copy column from table={\ThreeDRadialTwoCGLyapFalse}{[index] 2}] {timeCGLYAPinexact} {\ThreeDRadialTwoCG}
\pgfplotstablecreatecol[copy column from table={\ThreeDRadialTwoCGLyapFalse}{[index] 1}] {iterCGLYAPinexact} {\ThreeDRadialTwoCG}

\pgfplotstablecreatecol[copy column from table={\ThreeDRadialTwoGalerkin}{[index] 3}] {errorGalerkin} {\ThreeDRadialTwoGT}
\pgfplotstablecreatecol[copy column from table={\ThreeDRadialTwoGalerkin}{[index] 2}] {timeGalerkin} {\ThreeDRadialTwoGT}
\pgfplotstablecreatecol[copy column from table={\ThreeDRadialTwoGalerkin}{[index] 1}] {iterGalerkin} {\ThreeDRadialTwoGT}

\pgfplotstablecreatecol[copy column from table={\ThreeDRadialTwoTimestepping}{[index] 2}] {errorTimestepping} {\ThreeDRadialTwoGT}
\pgfplotstablecreatecol[copy column from table={\ThreeDRadialTwoTimestepping}{[index] 1}] {timeTimestepping} {\ThreeDRadialTwoGT}

\hspace{-1cm}

% Table with the Galerkin and time stepping values
\pgfplotstabletypeset[
%Rename the columns
%columns/Refinements/.style={column name={Refinement}},
%columns/errorCGOPTexact/.style={column name={$L_2$ error}},
%columns/timeCGOPTexact/.style={column name={Time}},
columns/errorCGOPTinexact/.style={column name={$L_2$ error}, column type = {|r}, sci},
columns/timeCGOPTinexact/.style={column name={Wall time}, column type = {|r}, sci},
columns/iterCGOPTinexact/.style={column name={Iter.}, column type = {|r}, fixed},
%columns/errorCGLYAPexact/.style={column name={$L_2$ error}},
%columns/timeCGLYAPexact/.style={column name={Time}},
columns/errorCGLYAPinexact/.style={column name={$L_2$ error}, column type = {|r}, sci},
columns/timeCGLYAPinexact/.style={column name={Wall time [s]}, column type = {|r}, sci},
columns/iterCGLYAPinexact/.style={column name={Iter.}, column type = {|r}, fixed},
columns/errorGalerkin/.style={column name={$L_2$-error}, column type = {|r}, sci},
columns/timeGalerkin/.style={column name={Wall time [s]}, column type = {|r}, sci},
columns/iterGalerkin/.style={column name={Iter.}, column type = {|r}, fixed},
columns/errorTimestepping/.style={column name={$L_2$-error}, column type = {|r}, sci},
columns/timeTimestepping/.style={column name={Wall time},column type/.add={|r}{|}, sci},
columns/Time/.style={column type = {|r},fixed},
columns/Space/.style={column type = {|r},fixed},
precision=3,
skip rows between index={0}{1},
% Head row		
every head row/.style={
	before row={
		\hline
		  \multicolumn{2}{|c|}{Unknowns} 
		 %& \multicolumn{2}{c|}{CG opt (exact)}
		 & \multicolumn{3}{c|}{Galerkin}
		 & \multicolumn{2}{c|}{Time stepping}\\
	},
	after row=\hline,		
},
column type/.add={|}{},
every last row/.style={after row=\hline}	
]{\ThreeDRadialTwoGT}

\caption{Case 2:  $L_2$-error and CPU time for the Galerkin projection and the time-stepping method.}
\label{table:radial_example_2GT}
\end{table}
\renewcommand{\arraystretch}{1.0} 

%----------------------------------------------------------------------------------
\subsection{Conclusions}
%----------------------------------------------------------------------------------
Our theoretical results and numerical experience show that the proposed  very weak variational space-time method is significantly more accurate than the Crank-Nicolson scheme on problems with low regularity, at comparable runtimes.

%===============================================
\begin{appendix}
%===============================================

%===============================================
\section{Proof of Theorem \ref{Thm:A.2}}
\label{App:A}
%===============================================
We collect the proof of the well-posedness for the semi-variational setting in \S\ref{SubSec:SemiVar}. 

\begin{prop}\label{Prop:A1}
	Let $s\in\R^+$. The mapping $w\mapsto \langle\cdot, w\rangle$, $w\in H^{-s}$, where
	\begin{align}\label{eq:Duality}
		\langle\cdot,\cdot\rangle: H^s\times H^{-s}\to\R,
		&\quad
		\langle v,w\rangle := \sum_{n=1}^\infty v_n\, w_n
	\end{align}
	is an isometric isomorphism from $H^{-s}$ to $(H^s)'$, i.e., $(H^s)'\cong H^{-s}$.
\end{prop}
\begin{proof}
	First, note that $H^s$ is a Hilbert space with the inner product $(v,w)_s := \sum_{n=1}^\infty \lambda_n^{s}\, v_n\, w_n$ and $H^s\hookrightarrow H \hookrightarrow H^{-s}$ with continuous embeddings. Let $v\in H^s$, $w\in H^{-s}$, then by H\"older's inequality
	\begin{align*}
		\langle v,w\rangle
		&= \sum_{n=1}^\infty \lambda_n^{s/2}v_n\, \lambda_n^{s/2}w_n 
		\le \left( \sum_{n=1}^\infty \lambda_n^{s}\, v_n^2\right)^{1/2}
			 \left( \sum_{n=1}^\infty \lambda_n^{s}\, w_n^2\right)^{1/2}
			= \| v\|_s\, \| w\|_{-s}<\infty.
	\end{align*}
	Hence, $\langle\cdot, w\rangle\in (H^s)'$ and $\| \langle\cdot,w\rangle\|_{(H^s)'} = \sup_{v\in H^s} \frac{\langle v,w\rangle}{\| v\|_{s}} \le \| w\|_{-s}$. 
	On the other hand, given $w\in H^{-s}$, set $\tilde v_n:=\lambda^{-s} w_n$ and $\tilde v:=\sum_{n=1}^\infty \tilde v_n e_n$. Then, 
$\| \tilde v\|_s^2 
		= \sum_{n=1}^\infty \lambda_n^s\, (\lambda_n^{-s}w_n)^2
		= \sum_{n=1}^\infty \lambda_n^{-s} (w_n)^2
		= \| w\|_{-s}^2 < \infty$,
	i.e., $\tilde v\in H^{-s}$. Moreover
$\langle \tilde v,w\rangle 
		= \sum_{n=1}^\infty \tilde v_n\, w_n
		= \sum_{n=1}^\infty \lambda_n^{-s} (w_n)^2
		= \| w\|_{-s}^2 = \| \tilde v\|_s\, \| w\|_{-s}$. 
	If $w\not=0$, we get that
$\|\langle\cdot,w\rangle\|_{(H^s)'} 
		= \sup_{v\in H^s} \frac{\langle v,w\rangle}{\| v\|_{s}} 
		\ge \frac{\langle \tilde v,w\rangle}{\| \tilde v\|_{s}} 
		= \| w\|_{-s}$ 
	with equality for $w=0$. Hence, $\|\langle\cdot,w\rangle\|_{(H^s)'}  = \| w\|_{-s}$ for all $w\in H^{-s}$.
\end{proof}

Now we start by considering the following \emph{homogeneous} abstract second order initial value problem. 
Let $u_0 \in D(A)$ and $u_1 \in H$. The goal is to find a function $w \in C^2([0,T], H)$ such that $w(t) \in D(A)$ for $t\in [0,T]$ and satisfying 
	\begin{align}
	\label{eq:9.4.2}
        \ddot w (t) + Aw(t) &= 0, \qquad t\in (0,T), 
        \qquad&
	w(0) = u_0, \ \dot w(0) = u_1,% \label{eq:9.4.2b} 	
	\end{align}
where the spaces $u_0$ and $u_1$ reside in will be specified later. It is easily seen that
(a) $u_0 = e_n$, $u_1 = 0$ yields $w(t) = \cos (\sqrt{\lambda_n} t) e_n$  and (b) $u_0 = 0$ and $u_1 = e_n$ gives rise to $w(t) = \lambda_n^{-1/2}\sin (\sqrt{\lambda_n} t) e_n$.
\smallskip

We can now express the general solution of \eqref{eq:9.4.2} as a series of solutions of these special types and prove the following theorem.
\smallskip

%---------------------------------------------------------------------------------------------------
\begin{thm}[Homogeneous wave equation]\label{Thm:A.1}
Let $s\in\R_{\ge 0}$, $u_0 \in H^s$ and $u_1 \in H^{s-1}$. Then \eqref{eq:9.4.2} admits a unique solution 	$w\in\cC^s$, see \eqref{eq:cCs}.
\end{thm}
%---------------------------------------------------------------------------------------------------
\begin{proof}
Uniqueness: Let $w\in\cC^s$ be a solution of \eqref{eq:9.4.2}, then $w(t)\in H$ for all $t\in [0,T]$. Set $w_n(t) := \langle w(t), e_n\rangle = (w(t), e_n)_H$ for $n\in\N$ and $t\in [0,T]$. Since $w\in\cC^s$, in particular $\ddot{w}(t) \in H^{s-2}$, we get by $e_n\in D(A)=H^2$ the fact 
$w_n\in C^2([0,T])$ with derivative
$\ddot w_n(t) 
		=  \langle \ddot w(t), e_n\rangle   = - \langle Aw (t), e_n\rangle 	
		= -\sum_{k=1}^{\infty} \lambda_k (w(t), e_k)_H \langle e_k, e_n\rangle 	
		= -\lambda_n (w(t), e_n)_H = - \lambda_n w_n (t)$, 
$t\in (0,T)$, 
and initial values
$w_n(0) = (w(0), e_n)_H = (u_0, e_n)_H$, 
$\dot w_n (0) 
		= \langle \dot{w}(0), e_n\rangle = \langle u_1, e_n\rangle$. 
This is an initial value problem of a second order linear ode with the unique solution
\begin{align}\label{eq:9.4.3a} 
	w_n (t) = \cos (\sqrt{\lambda_n} t) (u_0, e_n)_H 
		+ \lambda_n^{-1/2} \sin (\sqrt{\lambda_n} t) \langle u_1,e_n\rangle,
\end{align}
which is easily verified. Since 
\begin{align}\label{eq:9.4.3b} 
	w(t) = \sum_{n=1}^{\infty} w_n (t) e_n
\end{align}
is the unique expansion of $w(t)$ in $H$ with respect to the orthonormal basis $\{ e_n: n\in\N\}$, the uniqueness statement has been proved.

Existence: 
We now define $w_n(t)$ by \eqref{eq:9.4.3a} and \eqref{eq:9.4.3b}. Then, for all $t\in [0,T]$, 
\begin{align*}
	\| w(t)\|_s^2
%	&= \sum_{n=1}^\infty \lambda_n^s\, |w_n(t)|^2 \\
	&\le 2 \sum_{n=1}^\infty \lambda_n^s\, |\cos(t\sqrt{\lambda_n})|^2\, |(u_0,e_n)_H|^2
		+ 2 \sum_{n=1}^\infty \lambda_n^{s-1}\, |\sin(t\sqrt{\lambda_n})|^2\, |\langle u_1,e_n\rangle|^2 \\
	&\le 2 \sum_{n=1}^\infty \lambda_n^s\, |(u_0,e_n)_H|^2
		+ 2 \sum_{n=1}^\infty \lambda_n^{s-1}\, |\langle u_1,e_n\rangle|^2 	
	= 2\, \| u_0\|_s^2 + 2\, \| u_1\|_{s-1}^2 < \infty,
\end{align*}
uniformly in $t\in [0,T]$, so that $w\in C([0,T]; H^s)$. Next 
\begin{align*}
	\| \dot{w}(t)\|_{s-1}^2
%	&= \sum_{n=1}^\infty \lambda_n^{s-1}\, |\dot{w}_n(t)|^2 
	&\le 2 \sum_{n=1}^\infty \lambda_n^{s-1}\, \lambda_n |\sin(t\sqrt{\lambda_n})|^2\, |(u_0,e_n)_H|^2\\
	&\qquad\qquad%\qquad\qquad\qquad
		+ 2 \sum_{n=1}^\infty \lambda_n^{s-1}\, \lambda_n^{-1} \lambda_n |\cos(t\sqrt{\lambda_n})|^2\, |\langle u_1,e_n\rangle|^2 \\
	&\kern-8pt\le 2 \sum_{n=1}^\infty \lambda_n^s\, |(u_0,e_n)_H|^2
		+ 2 \sum_{n=1}^\infty \lambda_n^{s-1}\, |\langle u_1,e_n\rangle|^2 	
	= 2\, \| u_0\|_s^2 + 2\, \| u_1\|_{s-1}^2 < \infty,
\end{align*}
so that $w\in C^1([0,T]; H^{s-1})$ and similarly
\begin{align*}
	\| \ddot{w}(t)\|_{s-2}^2
%	&= \sum_{n=1}^\infty \lambda_n^{s-2}\, |\ddot{w}_n(t)|^2 
	&\le 2 \sum_{n=1}^\infty \lambda_n^{s-2}\, \lambda_n^2 |\cos(t\sqrt{\lambda_n})|^2\, |(u_0,e_n)_H|^2\\
	&\qquad\qquad%\qquad\qquad\qquad
		+ 2 \sum_{n=1}^\infty \lambda_n^{s-2}\, \lambda_n^{-1} \lambda_n^2 |\sin(t\sqrt{\lambda_n})|^2\, |\langle u_1,e_n\rangle|^2 \\
	&\kern-8pt\le 2 \sum_{n=1}^\infty \lambda_n^s\, |(u_0,e_n)_H|^2
		+ 2 \sum_{n=1}^\infty \lambda_n^{s-1}\, |\langle u_1,e_n\rangle|^2 	
	= 2\, \| u_0\|_s^2 + 2\, \| u_1\|_{s-1}^2 < \infty,
\end{align*}
which shows that $w\in C^2([0,T]; H^{s-2})$. We conclude that $w\in\cC^s$. Finally, we have
$\ddot w (t) 
 = \sum_{n=1}^{\infty} \ddot w_n (t) e_n = \sum_{n=1}^{\infty} w_n (t) \lambda_n e_n
 = -A w(t)$ 
by definition of $A$. In addition, 
$w (0) = \sum_{n=1}^{\infty} (u_0, e_n)_H e_n = u_0$ and
$\dot w (0)  = \sum_{n=1}^{\infty} \langle u_1, e_n\rangle\, e_n = u_1$.
This shows that $w$ solves \eqref{eq:9.4.2}, and we have proved existence of solutions.
\end{proof}

We are now in the position to prove Theorem~\ref{Thm:A.2} for the wave equation with inhomogeneous right-hand side.

\begin{proof}[Proof of Theorem \ref{Thm:A.2}.]
	Since the difference of two solutions of \eqref{Eq:8.37neu} is a solution of the homogeneous problem \eqref{eq:9.4.2}, uniqueness follows from Theorem \ref{Thm:A.1}. Moreover, since the homogeneous problem has a solution, in order to prove existence for \eqref{Eq:8.37neu}, we may and will assume that $u_0=u_1=0$. 
	
	Next, we set $f_n(t):= \langle f(t), e_n\rangle $, which is well-defined since $e_n\in D(A)=H^2$ and $f(t)\in H^{s-1}$, $s\ge 0$. Then, $f_n\in C([0,T])$. We set 
$w_n(t) = \lambda_n^{-1/2} \int_0^t \sin(\sqrt{\lambda_n} (t-\tau))\, f_n(\tau)\, d\tau$, 
	and $w(t):= \sum\limits_{n=1}^\infty \langle w_n(t), e_n\rangle\, e_n$.  By H\"older's inequality, we have, for all $t\in [0,T]$ 
	\begin{align*}
		\lambda_n^s\, w_n(t)^2
		&\le \lambda_n^{s-1} \int_0^T \sin(\sqrt{\lambda_n} (t-\tau))^2\, d\tau \int_0^T f_n(\tau)^2\, d\tau
		\le T \lambda_n^{s-1} \, \int_0^T f_n(\tau)^2\, d\tau,
	\end{align*} 
	so that
	\begin{align*}
		\| w(t)\|_s^2
%		&= \sum_{n=1}^\infty \lambda_n^s\, w_n(t)^2\\
		&\le T \sum_{n=1}^\infty \lambda_n^{s-1} \kern-5pt\int_0^T \kern-4pt f_n(\tau)^2\, d\tau
		= T\int_0^T \sum_{n=1}^\infty \lambda_n^{s-1} f_n(\tau)^2\, d\tau%\\
		= T\int_0^T \kern-3pt\| f(\tau)\|_{s-1}^2\, d\tau,
	\end{align*}
	which is finite uniformly in $t\in [0,T]$ since $f\in C([0,T]; H^{s-1})$, so that $w\in C([0,T]; H^s)$. Next, we note that 
$\dot{w}_n(t) = \int_0^t \cos( \sqrt{\lambda_n} (t-s))\, f_n(s)\, ds$, 
	so that similar as above
	\begin{align*}
		\lambda_n^{s-1}\, w_n(t)^2
		&\le \lambda_n^{s-1} \int_0^T \kern-4pt\cos(\sqrt{\lambda_n} (t-\tau))^2\, d\tau \int_0^T \kern-4pt f_n(\tau)^2\, d\tau
		\le T \lambda_n^{s-1} \, \int_0^T \kern-4pt  f_n(\tau)^2\, d\tau,
	\end{align*} 
	which yields 
	\begin{align*}
		\| \dot{w}(t)\|_{s-1}^2
%		&= \sum_{n=1}^\infty \lambda_n^{s-1}\, \dot{w}_n(t)^2\\
		&\kern-2pt\le T \sum_{n=1}^\infty \lambda_n^{s-1}\kern-4pt \int_0^T\kern-4pt  f_n(\tau)^2\, d\tau
		= T\kern-2pt\int_0^T \kern-2pt\sum_{n=1}^\infty \lambda_n^{s-1} f_n(\tau)^2\, d\tau
		= T\kern-2pt\int_0^T\kern-4pt  \| f(\tau)\|_{s-1}^2\, d\tau,
	\end{align*}
	which again is finite uniformly in $t\in [0,T]$, so that $w\in C^1([0,T]; H^{s-1})$. In order to prove $w\in C^2([0,T]; H^{s-2})$ (and thus $w\in\cC^s$), we note that $\ddot{w}_n+\lambda_n\, w_n = f_n$, $w_n(0)=\dot{w}_n(0)=0$. Hence, 
	\begin{align*}
		\| \ddot{w}(t)\|_{s-2}^2
%		&= \sum_{n=1}^\infty \lambda_n^{s-2}\, \ddot{w}_n(t)^2\\
		&\le 2 \sum_{n=1}^\infty \lambda_n^{s-2}\lambda_n^2 w_n(t)^2
			+ 2 \sum_{n=1}^\infty \lambda_n^{s-2} f_n(t)^2
		= 2\, \| w(t)\|_s^2 + 2\, \| f(t)\|_{s-2}^2 < \infty 
	\end{align*}
	uniformly in $t\in [0,T]$, so that $w\in C^2([0,T]; H^{s-2})$. Finally
	\begin{align*}
		\ddot{w}(t)
		&= \sum_{n=1}^\infty \langle \ddot{w}_n(t), e_n\rangle\, e_n
		= - \sum_{n=1}^\infty \lambda_n\, \langle w_n(t), e_n\rangle\, e_n
			+ \sum_{n=1}^\infty \langle f_n(t), e_n\rangle\, e_n \\
		&=- \sum_{n=1}^\infty \lambda_n\, (w_n(t), e_n)_H\, e_n + f(t)
		= -Aw(t) + f(t)
	\end{align*}
	for all $t\in (0,T)$. Since $w_n(0) = 0 = \dot{w}_n(0)$, we obtain $w(0)=\dot{w}(0)=0$, so that $w$ solves \eqref{Eq:8.37neu} for $u_0=u_1=0$, which concludes the proof.
\end{proof}

\end{appendix}

%%========================================================================
\bibliographystyle{abbrv}
\bibliography{ST_references}
%%========================================================================
\end{document}